\documentclass[11pt,a4paper]{amsart}
\usepackage{enumitem}
\usepackage[english]{babel}
\usepackage{tikz-cd}
\usepackage[margin=2.7cm]{geometry}
\usepackage{array}
\usepackage{multirow}
\usepackage{amsfonts, amssymb, amsmath, amsthm, stmaryrd, latexsym}
\usepackage{graphicx}
\usepackage{epstopdf}
\usepackage{algorithmic}
\usepackage{relsize}
\usepackage{amsopn}
\usepackage{mathrsfs}
\usepackage[cmtip,arrow]{xy}
\usepackage{pb-diagram,pb-xy}
\usepackage{amscd,extarrows}
\usepackage{tikz}
\usepackage{subcaption}
\usepackage{csquotes}
\usepackage{graphicx}
\usepackage{bm}
\usepackage[labelformat=simple, labelfont=normal]{subfig}

\theoremstyle{plain}
\newtheorem{thm}{Theorem}[section]

\theoremstyle{definition}
\newtheorem{defn}[thm]{Definition}

\theoremstyle{remark}
\newtheorem{remark}[thm]{Remark}
\newtheorem{example}[thm]{Example}

\usetikzlibrary{decorations.markings}
\tikzstyle arrowstyle=[scale=2]
\tikzstyle directed=[postaction={decorate,decoration={markings,
    mark=at position 0.55 with {\arrow[arrowstyle]{stealth}}}}]
\usepackage[
    backend=biber,
    style=numeric,
  ]{biblatex}

\ifpdf
  \DeclareGraphicsExtensions{.eps,.pdf,.png,.jpg}
\else
  \DeclareGraphicsExtensions{.eps}
\fi

\usepackage{enumitem}
\setlist[enumerate]{leftmargin=.5in}
\setlist[itemize]{leftmargin=.5in}

\title{Data-driven Identification of Attractors using Machine Learning}

\author[Gameiro]{Marcio Gameiro}
\address{Department of Mathematics, Rutgers University}
\email{gameiro@math.rutgers.edu}

\author[Gelb]{Brittany Gelb}
\address{Department of Mathematics, Rutgers University}
\email{brittany.gelb@rutgers.edu}

\author[Kalies]{William D.~Kalies}
\address{Department of Mathematics and Statistics, University of Toledo}
\email{william.kalies@utoledo.edu}

\author[Kram\'ar]{Miroslav Kram\'ar}
\address{Department of Mathematics, University of Oklahoma}
\email{miro@ou.edu}

\author[Mischaikow]{Konstantin~Mischaikow}
\address{Department of Mathematics, Rutgers University}
\email{mischaik@math.rutgers.edu}

\author[Tatasciore]{Paul Tatasciore}
\address{Department of Physics, Florida Atlantic University}
\email{ptatasciore2016@fau.edu}

\DeclareMathOperator{\cl}{cl}
\DeclareMathOperator{\Int}{int}

\newcommand{\sA}{{\mathsf A}}

\newcommand{\sN}{{\mathsf N}}

\newcommand{\R}{{\mathbb{R}}}
\newcommand{\Z}{{\mathbb{Z}}}

\newcommand{\cI}{{\mathcal I}}

\begin{document}

\begin{abstract}
In this paper we explore challenges in developing a topological framework in which machine learning can be used to robustly characterize global dynamics.
Specifically, we focus on learning a useful discretization of the phase space of a flow on compact, hyperrectangle in $\R^n$ from a neural network trained on labeled orbit data.
A characterization of the structure of the global dynamics is obtained
from approximations of attracting neighborhoods provided by the phase space discretization. The perspective that motivates this work is based on Conley's topological approach to dynamics, which provides a means to evaluate the efficacy and efficiency of our approach. 
\end{abstract}

\maketitle

\section{Introduction}

The analysis of applied dynamical systems is undergoing a paradigm shift toward data-driven methods. This change has been triggered by an explosion of data available for large, multiscale problems, as well as the development of powerful, computational methods for data-based analysis 
facilitated by machine learning.
Nevertheless, the fundamental questions remain the same: how do we achieve an insightful characterization of the dynamics, and how can we predict the dynamics?

Current methods often fit the data to an explicit model that takes the form of a parameterized differential equation or difference equation \cite{SINDY,Ives08}.
The clear advantage of this approach is that it allows the user to exploit a wide variety of well-developed techniques to analyze the associated dynamics. 
Alternatively, black-box surrogate models can  be learned directly from the data.
This approach typically provides for less interpretability of the model, but may lead to efficient predictability of the dynamics \cite{ren2021deep, schultz2021can}. However, bifurcation theory suggests that descriptions of long-term dynamics provided by either of these approaches may not be robust.

The perspective that motivates this work is based on Conley's topological approach to dynamics \cite{Conley}, which provides a framework for deducing global information about dynamics that is robust with respect to perturbations.
This framework has been successfully applied to the analysis of time series data \cite{batko, mischaikow:mrozek:reiss:szymcak}, machine-learned dynamics \cite{vieira123}, and rigorous verification of numerical analysis of differential equations \cite{vandenberg:lessard, gameiro2024globaldynODE, mischaikow:mrozek:95} and difference equations \cite{Database_schema, day:junge:mischaikow}.
Two steps are common to most of these applications: ({\bf Step 1}) identification of regions in which the dynamics of interest occurs, and ({\bf Step 2}) use of the Conley index, an algebraic topological invariant, to identify the existence and structure of the dynamics in the regions of interest. 
Theoretical results \cite{gameiro:gelb:mischaikow} indicate that with sufficient data and sufficient computational resources, both {\bf Step 1} and {\bf Step 2} can be accessed via machine learning.
However, \cite{gameiro:gelb:mischaikow} provides no quantification as to how much data or how much computation is required.

In this paper we undertake an initial exploration of the challenges of using this topological framework to machine learn dynamics and focus on {\bf Step 1}, identifying regions of interest.
As is discussed in Section~\ref{sec:dynamics}, all regions of interest can be derived from a lattice of attracting neighborhoods.
With this in mind, we work with systems that have multiple global basins of attraction and explore the use of machine learning to identify discretizations of phase space that provide approximations of attracting neighborhoods from which the correct Conley index information can be computed. 

Using various examples, with dynamics coming from known ordinary differential equations (ODEs), we show that if test loss is sufficiently low, then the Conley indices of the attracting neighborhoods match the Conley indices of the known attractors. 
This is significant because, as is mentioned above, the information about the structure of dynamics is deduced using the Conley index. 
In this sense we view this work as a proof of concept that machine learning and computational Conley theory can be used together for applications in nonlinear dynamics.\footnote{Code is available at: https://github.com/begelb/MLCD}

However, there are important caveats to this claim. To compute the Conley index we need to be able to compute homology. Computational homology is a relatively new endeavor \cite{kaczynski:mischaikow:mrozek}, and while considerable progress is being made \cite{fugacci2014efficient, harker:mischaikow:mrozek:nanda},  to the best of our knowledge the most manageable package for dealing with complexes of moderate dimension is \texttt{pyCHomP} \cite{pyCHomP_repo}. Since \texttt{pyChomP} is most efficient on cubical complexes, we restrict our study to special neural networks to achieve a decomposition of the space that fits into the data structure of a cubical grid. We use the strategy mentioned in~\cite{Balestriero2019} to ensure that the network partitions the space into parallelotopes that can be transformed to cubes by a linear coordinate change.  There is a price to be paid for using these networks, namely they cannot efficiently approximate attracting neighborhoods with complicated geometry.  However, we  demonstrate on several examples that as long as the network is a sufficiently good approximation to the labeling function introduced in Section~\ref{sec:AI}, then we obtain the correct Conley indices. 

We also present a few examples where the geometry is too complicated to be captured properly by the restricted network. For these examples we present heuristics that detect if the network is versatile enough to capture the geometry of the basins of attraction. Our heuristics  apply to unrestrained networks as well, which is critical for further 
development of the proposed framework. Finally, we show that the number of grid elements (parallelotopes) in the machine learned decomposition is radically smaller than in a regular grid. This observation is crucial because despite the fast improvements in computational homology, a large number of cubes is still a serious bottleneck. 

The rest of the paper is organized as follows. In Section~\ref{sec:dynamics}, we recall the basic concepts of the Conley-Morse theory for ordinary differential equations. Development of our method in Section~\ref{sec:DDidentification} is divided into two parts. In Section~\ref{sec:AI}, we present possible strategies for identifying attractors and creating a sample of a labeling function that describes to which basin of attraction a given point belongs. In Section~\ref{sec:ML}, we use a neural network to  partition the domain of the ODE into polytopes. We also explain there how to use these polytopes to approximate the attracting neighborhoods of the attractors identified in Section~\ref{sec:AI}. These attracting neighborhoods are crucial for computing the Conley indices. The details about training the neural networks and benchmarking our method are presented in Section~\ref{sec:training}. The results are divided into three parts. Section~\ref{sec:fixed_points} familiarizes the reader with the pipeline of our method on simple examples, and Section~\ref{sec:periodic_orbits} highlights the reduction in the size of the grid produced by our method. Two systems that exceed the versatility of the restrained network are discussed in Section~\ref{sec:complex}, where we also provide heuristics on deciding if the network is sufficiently expressive to produce the correct Conley indices. Finally in Section~\ref{sec:conclusion}, we outline how our framework can be further developed to overcome the limitations of this initial study.

\section{Dynamics and Conley Index Theory}
\label{sec:dynamics}

We are interested in characterizing the dynamics of an ODE 
\begin{equation}
    \label{eq:basicODE}
    \frac{dx}{dt} = g(x)
\end{equation} 
on a compact region $X\subset \R^d$.
For the purposes of this paper we choose $X$ to be a hyperrectangle of the form $\prod_{i = 1}^d [a_i, b_i]$.
For simplicity of exposition we assume that the solutions to the ODE give rise to a \emph{flow} $\varphi\colon \R \times \R^d \to \R^d$; that is, a continuous function satisfying $\varphi(0,x) = x$  and $\varphi(t, \varphi(s, x)) = \varphi(t + s,x)$ for all $t, s \in \R$ and $x \in \R^d$.
Our description of the dynamics of \eqref{eq:basicODE} involves identifying \emph{invariant sets}; that is, sets $S\subset X$ such that $\varphi(t,S)=S$ for all $t\in\R$.

As indicated in the introduction, {\bf Step 1} involves the identification of regions of phase space in which dynamics of interest occurs.
Ab initio these regions are not known and thus must be identified via the dynamics.
Perhaps the most obvious region to attempt to identify is an \emph{attracting neighborhood}; that is a compact set $N\subset X$ with the property that there exists $\tau > 0$ such that $\varphi(t,N)\subset \Int(N)$ for all $t \geq \tau$, where $\Int$ denotes interior.
Observe that the condition that $N$ is compact and $\varphi(\tau,N)\subset \Int(N)$ guarantees that an attracting neighborhood is \emph{robust}; that is, if $N$ is an attracting neighborhood for $\frac{dx}{dt} = g(x)$, then $N$ remains an attracting neighborhood even if the nonlinearity $g$ is slightly perturbed.

Attracting neighborhoods give rise to invariant sets.
In particular, if $N$ is an attracting neighborhood, then
\[
A = \omega(N) := \bigcap_{t > 0}\cl(\varphi([t, \infty),N))
\]
is an invariant set. 
An invariant set $A$ of this form is called an \emph{attractor}, and $N$ is referred to as an attracting neighborhood of $A$.
The \emph{basin of attraction} of an attractor consists of the union of all of its attracting neighborhoods.
This last remark emphasizes an important point; an attracting neighborhood $N$ identifies a unique attractor $A =\omega(N)$, but an attractor typically has uncountably many attracting neighborhoods.
Thus, at least abstractly, it is easier to characterize dynamics in terms of attractors than in terms of attracting neighborhoods.

Observe that identification of an attracting neighborhood $N$ provides  predictive power; an initial condition in $N$ asymptotically approaches the attractor $A = \omega(N)$. 
In particular, we leave it to the reader to check that if $A$ is an attractor and $x\in X$, then $\omega(x) \subset A$ if and only if $x$ is an element of the basin of attraction of $A$.
The dynamics on an attractor can be quite complicated, for example chaotic dynamics on strange attractors. More importantly, bifurcation theory tells us that the dynamics can change dramatically even with respect to small perturbations of the system. 
This implies that direct numerical simulation cannot in general provide a precise characterization of the dynamics on an attractor. 
However, as indicated above, attracting neighborhoods are robust with respect to perturbations, and thus together with Conley index information can provide coarse but accurate information about dynamics within the attracting neighborhood.

As indicated in the introduction, we focus on multistable systems, and thus we need a language in which to organize the structure of attractors. As is discussed in \cite{LSoA1}, the collection of attractors of a dynamical system has the structure of a bounded, distributive lattice with partial order given by inclusion.
The attractor lattice of a system need not be finite, but for systems in $\R^d$, it is countable \cite{robinson}. 
Furthermore, finite sublattices of attractors can be used to approximate the dynamics at finer and finer scales \cite{LSoA2,LSoA3}. In this setting, a system exhibits \textit{multistability} if there is a finite sublattice $\mathsf{A}$ of attractors that contains more than one \textit{minimal} attractor, i.e.\ an attractor with no nonempty predecessor in $\mathsf{A}$.

By \cite[Proposition~4.3]{LSoA1}, the collection of attracting neighborhoods forms a bounded distributive lattice with partial order given by inclusion, and $\omega$ maps lattices of attracting neighborhoods to lattices of attractors homomorphically. From this structure, one can deduce the following result that constructs $\mathsf{A}$ from multiple, distinct attracting neighborhoods.

\begin{thm}
\label{thm:fundamental}
Suppose that $X$ is an attracting neighborhood for \eqref{eq:basicODE}.
Let $\{N_k \subset X\mid k = 0,\ldots, K\}$ be a set of mutually disjoint attracting neighborhoods with $A_k := \omega(N_k)$.
Then
\[
\sA := \omega(X) \cup \left\{ \bigcup_{k\in I}A_k \mid I\subset \{0,\ldots, K \} \right\} 
\]
is a finite lattice of attractors with minimal attractors $A_k$.
\end{thm}

The importance of Theorem~\ref{thm:fundamental} for this work is that from data we machine learn a collection of disjoint attracting neighborhoods, and therefore identify a lattice of attractors.
By definition of the lattice, the number of elements in $\sA$ grows rapidly with the number of disjoint attracting neighborhoods $K$. Thus, rather than working with $\sA$ we turn to an equivalent representation that takes the form of a partially ordered set (poset).

\begin{defn}
\label{def:Morse}
{\em
A \emph{Morse representation} of a dynamical system on a compact metric space $X$ is a finite, partially ordered set $\{M_i\}_{i=1,\ldots n}$ of nonempty, pairwise-disjoint, compact, isolated invariant sets such that for all $x \in \omega(X) ~\backslash~ \bigcup\limits_{i = 1}^{n} M_i$ there exist $M_i <M_j$ such that $\omega (x) \subset M_i$ and $\alpha (x) \subset M_j$, where $\alpha(x):= \bigcap_{t>0}\cl(\varphi(x, (-\infty, -t]))$. 
The elements of a Morse representation are called \emph{Morse sets}.}
\end{defn}

Even a cursory explanation of the general correspondence between lattices of attractors and  Morse representations is beyond the scope of this article (the interested  reader is referred to  \cite{LSoA3}).
Instead, we consider a specific example that we hope provides sufficient intuition. 

\vskip 12pt

\begin{figure*}[!htb]
\centering
\begin{subfigure}{0.36\textwidth}
\centering
\begin{tikzpicture}
[main node/.style={
},
scale=1.55]
\node[main node] (1) at (0,0) {$\emptyset$};
\node[main node] (2) at (-1,0.5) {$N_0$};
\node[main node] (3) at (0,0.5) {$N_1$};
\node[main node] (4) at (1,0.5) {$N_2$};
\node[main node] (5) at (-1,1) {$N_0\cup N_1$};
\node[main node] (6) at (0,1) {$N_0\cup N_2$};
\node[main node] (7) at (1,1) {$N_1\cup N_2$};
\node[main node] (8) at (0,1.5) {$N_0\cup N_1\cup N_2$};
\node[main node] (9) at (0,2) {$X$};
\path[thick]
(1) edge[-, shorten <= 2pt, shorten >= 2pt] (2)
(1) edge[-, shorten <= 2pt, shorten >= 2pt] (3)
(1) edge[-, shorten <= 2pt, shorten >= 2pt] (4)
(2) edge[-, shorten <= 2pt, shorten >= 2pt] (5)
(2) edge[-, shorten <= 2pt, shorten >= 2pt] (6)
(3) edge[-, shorten <= 2pt, shorten >= 2pt] (5)
(3) edge[-, shorten <= 2pt, shorten >= 2pt] (7)
(4) edge[-, shorten <= 2pt, shorten >= 2pt] (6)
(4) edge[-, shorten <= 2pt, shorten >= 2pt] (7)
(5) edge[-, shorten <= 2pt, shorten >= 2pt] (8)
(6) edge[-, shorten <= 2pt, shorten >= 2pt] (8)
(7) edge[-, shorten <= 2pt, shorten >= 2pt] (8) 
(8) edge[-, shorten <= 2pt, shorten >= 2pt] (9);
\end{tikzpicture}
\caption{}
\end{subfigure}
\centering
\begin{subfigure}{0.36\textwidth}
\centering
\begin{tikzpicture}
[main node/.style={
},
scale=1.55]
\node[main node] (1) at (0,0) {$\emptyset$};
\node[main node] (2) at (-1,0.5) {$A_0$};
\node[main node] (3) at (0,0.5) {$A_1$};
\node[main node] (4) at (1,0.5) {$A_2$};
\node[main node] (5) at (-1,1) {$A_0\cup A_1$};
\node[main node] (6) at (0,1) {$A_0\cup A_2$};
\node[main node] (7) at (1,1) {$A_1\cup A_2$};
\node[main node] (8) at (0,1.5) {$A_0\cup A_1\cup A_2$};
\node[main node] (9) at (0,2) {$\omega(X)$};
\path[thick]
(1) edge[-, shorten <= 2pt, shorten >= 2pt] (2)
(1) edge[-, shorten <= 2pt, shorten >= 2pt] (3)
(1) edge[-, shorten <= 2pt, shorten >= 2pt] (4)
(2) edge[-, shorten <= 2pt, shorten >= 2pt] (5)
(2) edge[-, shorten <= 2pt, shorten >= 2pt] (6)
(3) edge[-, shorten <= 2pt, shorten >= 2pt] (5)
(3) edge[-, shorten <= 2pt, shorten >= 2pt] (7)
(4) edge[-, shorten <= 2pt, shorten >= 2pt] (6)
(4) edge[-, shorten <= 2pt, shorten >= 2pt] (7)
(5) edge[-, shorten <= 2pt, shorten >= 2pt] (8)
(6) edge[-, shorten <= 2pt, shorten >= 2pt] (8)
(7) edge[-, shorten <= 2pt, shorten >= 2pt] (8) 
(8) edge[-, shorten <= 2pt, shorten >= 2pt] (9);
\end{tikzpicture}
\caption{}
\end{subfigure}
\begin{subfigure}{0.25\textwidth}
\centering
\begin{tikzpicture}
[main node/.style={
},
scale=2.0]
\node[main node] (2) at (-0.5,0.5) {$M_0$};
\node[main node] (3) at (0,0.5) {$M_1$};
\node[main node] (4) at (0.5,0.5) {$M_2$};
\node[main node] (6) at (0,1) {$M$};
\path[thick]
(2) edge[-, shorten <= 2pt, shorten >= 2pt] (6)
(3) edge[-, shorten <= 2pt, shorten >= 2pt] (6)
(4) edge[-, shorten <= 2pt, shorten >= 2pt] (6);
\end{tikzpicture}
\caption{}
\end{subfigure}
\caption{(a) Lattice of attracting neighborhoods derived from three mutually disjoint attracting neighborhoods. (b) Associated lattice of attractors. (c) Associated Morse representation with partial order relations $M_k < M$ and $M_k$ not related to $M_i$ for $i,k\in \{0,1,2 \}$.}
\label{fig:attMorse}
\end{figure*}
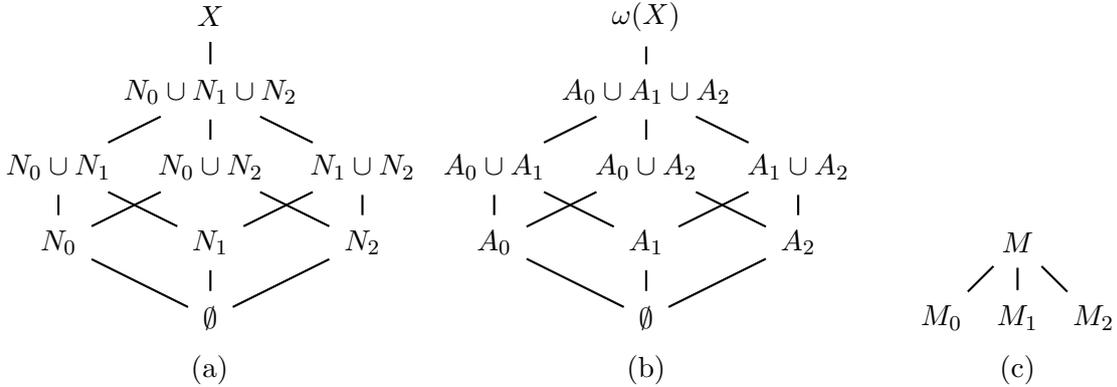

\begin{example}
\label{ex:attMorse}
{\em
Consider a hyperrectangle $X\subset \R^d$ and assume it is an attracting neighborhood for a flow $\varphi\colon \R\times \R^d \to \R^d$.
Let $\{N_k \subset X \mid k=0,1,2\}$ be mutually disjoint attracting neighborhoods.
As discussed in Section~\ref{sec:DDidentification}, we use machine learning to identify these sets.
Figure~\ref{fig:attMorse}(a) shows the associated lattice of attracting neighborhoods.
Let $A_k := \omega(N_k)$.
A topological argument ($X$ is connected and hence $\omega(X)$ is connected) implies that $\omega(X) \neq A_0\cup A_1\cup A_2$.
By Theorem~\ref{thm:fundamental}, we obtain the lattice of attractors $\sA$ shown in Figure~\ref{fig:attMorse}(b).

The associated Morse representation is shown in Figure~\ref{fig:attMorse}(c).
The minimal Morse sets $\{M_k\mid k=0,1,2 \}$ are the attractors, i.e., $M_k = A_k$.
The remaining Morse set is defined as follows
\[
M := \{ x\in \omega(X) \mid \omega(x)\cap M_k = \emptyset\ \text{for all $k=0,1,2$}\}.
\]
A point set topological argument based on the continuity of the flow $\varphi$ shows that $M\neq \emptyset$, $M$ is compact, and $M$ is disjoint from $M_k$ for $k=0,1,2$. To justify the partial order relation, we note that if $x\in N_k$, then $\omega(x)\subset A_k$.
Therefore, if $x\in \omega(X),$ then it is  impossible for $\alpha(x)\subset M_k$ and $\omega(x)\subset M$ for any $k=0,1,2$.}
\end{example}

Observe that in Example~\ref{ex:attMorse}, we explicitly used topology to argue that $\omega(X) \neq A_0\cup A_1\cup A_2$ or equivalently that $M\neq \emptyset$.
More sophisticated topological arguments allow us to extract additional information concerning the structure of the dynamics.
In particular, the \textit{Conley index}, a homological invariant, can be used to prove the existence of invariant sets and provide information about their structure.
The most fundamental result is that  if a Conley index  is nontrivial, then there is an associated nonempty invariant set  \cite{Conley}.
However, it can also be used to prove the existence of equilibria \cite{srzednicki85,mccord:88}, periodic orbits \cite{mccord:mischaikow:mrozek}, heteroclinic orbits \cite{Conley}, and chaotic dynamics \cite{mischaikow:mrozek:95,szymczak96,  day:frongillo}. 

The Conley index is defined to be the relative homology of an index pair \cite{Conley}.
Rather than providing the most general definition, we restrict our attention to a special case that is related to the computational efforts of this paper.
An \emph{attracting block} is an attracting neighborhood $N$ with the property that $\varphi(t,N)\subset \Int(N)$ for all $t >0$.
As in the case of attracting neighborhoods, attracting blocks form a bounded distributive lattice with partial order given by inclusion.

\begin{remark}
\label{rem:index}
Assume that we are given a finite lattice of attracting blocks $\sN$.
If $N_0,N_1 \in \sN$ and $N_0\subset N_1$, then $H_*(N_1,N_0)$ is the Conley index for the maximal invariant set in $N_1\setminus N_0$. In the case that $N_1$ is minimal, then $N_0=\emptyset$, as in Figure~\ref{fig:attMorse}(a). Thus the Conley index is simply given by $H_*(N_1).$ 

For the purpose of this paper we compute homology using $\Z_2$ coefficients.
Therefore, the Conley index is a sequence of vector spaces, and so if $X\subset \R^d$, it is sufficient to express it via the sequence of the Betti numbers, i.e.\ the dimensions of the vector spaces, $(\beta_0,\beta_1, \ldots, \beta_d)$.
\end{remark}

\section{Data-Driven Identification of Attracting Neighborhoods}
\label{sec:DDidentification}

The goal of this section is to describe a data-driven computational scheme for identifying attracting neighborhoods.\footnote{Data used for the paper is available at: https://github.com/begelb/MLCD-data} The input is a collection of orbit segments from a finite set of initial conditions $\cI \subset X$. In Section~\ref{sec:AI}, we use the orbit segments to identify finite approximations of attractors $A_k$ for $k\in\{0,\ldots,K\}$. This produces a labeling function $F \colon \cI \to \{0,\ldots, K\}$ that maps each initial condition to an attractor. In Section~\ref{sec:ML}, we discuss the use of a neural network to extend the labeling function $F$ to all of $X$ and how the trained network can be used to identify an attracting neighborhood $N_k$ for each $A_k$. In Section~\ref{sec:training} we explain the neural network training.

\subsection{Attractor identification scheme}
\label{sec:AI}

\begin{figure}
\centering
\begin{subfigure}{.33\textwidth}
  \centering
  \includegraphics[width=.8\linewidth]{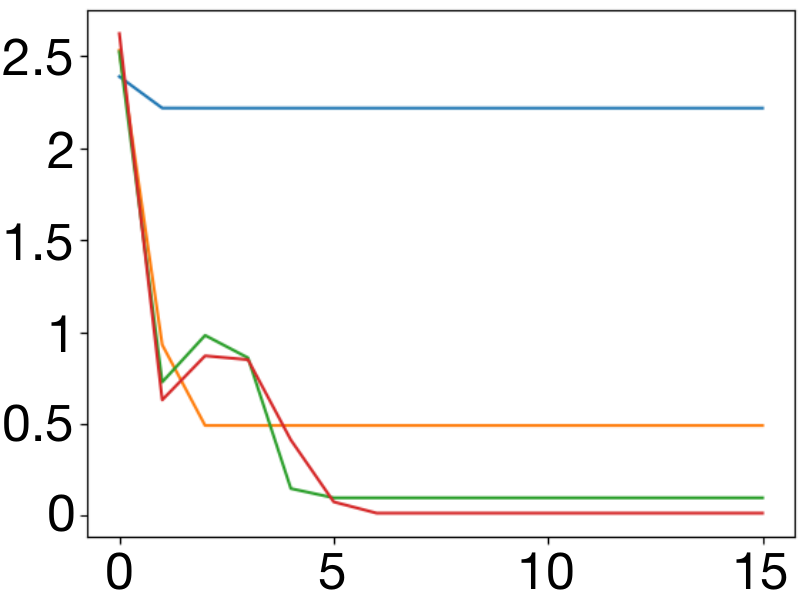}
  \caption{}
\end{subfigure}%
\begin{subfigure}{.33\textwidth}
  \centering
  \includegraphics[width=.6\linewidth]{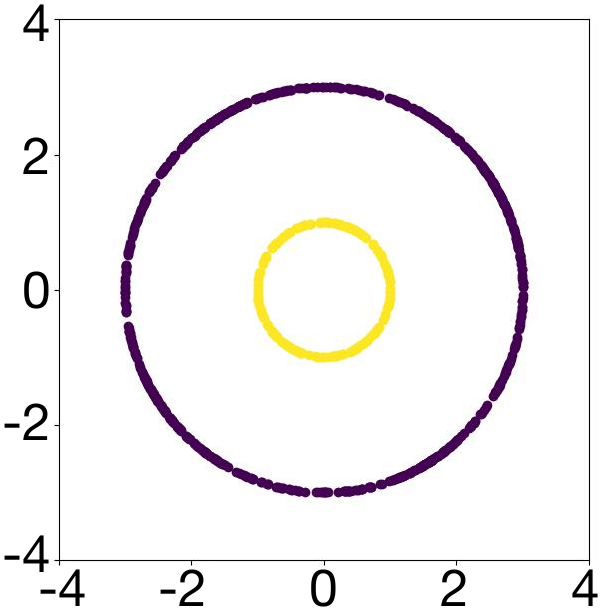}
  \caption{}
\end{subfigure}
\begin{subfigure}{.33\textwidth}
  \centering
  \includegraphics[width=.6
  \linewidth]{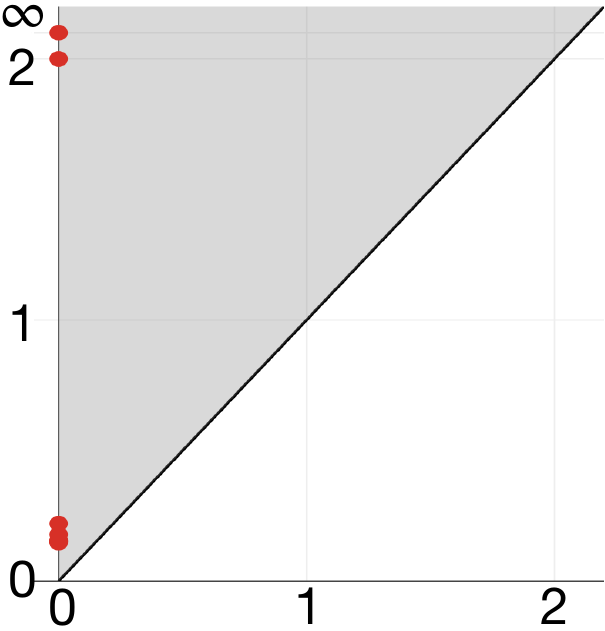}
  \caption{}
\end{subfigure}
 \caption{(a) Hausdorff distance between iterations for $10$ points (blue), $100$ points (orange), $1,000$ points (green), and $10,000$ points (red) for the radially symmetric flow defined by Equation~\eqref{eqn:radial_2labels1}. (b) The image of  $f_1^{14}(\cI)$ for $\cI$ consisting of $1,000$ points. (c) The zero-dimensional persistence diagram of the Vietoris-Rips complex with vertices given by $f_1^{14}(\cI)$, shown in part (b). Two clearly separated points which are far from the diagonal indicate the presence of two distinct clusters.
 }
\label{fig:hd_vary_pts}
\end{figure}

The input for the data-driven identification  of the attracting neighborhoods is a collection of orbit segments of a given ordinary differential equation (ODE). In this paper we consider the ODE on a hyperrectangle $X$, and the orbit segments start from initial conditions in $\cI \subset X$. We use  numerical simulation to approximate orbit segments of an explicit system, but in principle for physical systems, they could be obtained from experimental data.   The set $\cI$ of the initial conditions for the input orbit segments can be a random sample from $X$, set of vertices of some grid, or a set constructed by some space-filling design.  In this paper, we apply a Latin hypercube space-filling design~\cite{SFD} and perform the sampling by using the method {\tt LatinHypercube} implemented in the {\tt scipy.stats.qmc} module. 

Different choices of  $\cI$ might lead to identification of different sets $A_k$ and $N_k$, but Theorem~\ref{thm:fundamental} guarantees that the constructed lattice of attractors is correct. However, if the sample $\cI$ is too small then the attractors might not be correctly identified. Similarly, learning of the labeling function $F$ can be negatively affected by a small sample size. In this section we provide some heuristics to assess if the sample size is sufficient, but we do not conduct a systematic investigation of the sample size necessary to properly train the network. To avoid this issue, we use a reasonably large sample as explained in Section~\ref{sec:net_params}.

To generate samples of the orbit segments starting from initial conditions in the given set $\cI$, we use a standard numerical ODE solver {\tt odeint} implemented in the {\tt scipy.integrate} library. For each initial condition $x_0 \in \cI$ we construct a sample \{$x(i)\}_{i=0}^I$, of the orbit $x$ satisfying the ODE and initial value condition $x(0) = x_0$. The number of time steps $I$ is the same for all the orbits. The sampled orbits provide an  approximation of the time-$1$ map $f_1$ and its iterates $f_1^i$ for $x \in \cI$ and $0 \leq i \leq I$. To identify the asymptotic behavior from the sampled orbits we suppose that at time $I$ the orbits are close to their $\omega$-limit sets. To check if this assumption is reasonable, for our choice of the time horizon $I$, we consider the Hausdorff distance $d_H(f_1^i(\cI), f_1^{i+1}(\cI))$ between the images of the set $\cI$ under consecutive iterates of the time-$1$ map. In particular if the set $\cI$ is sufficiently dense and the time $I$ is long enough to ensure that the orbits are already close to their attractors, then $d_H(f_1^i(\cI),f_1^{i+1}(\cI))$ should be small. The appropriate size of $\cI$ and the value $I$ are system dependent. However, since the goal is to identify attracting neighborhoods rather than precise dynamics on the attractors, there is considerable flexibility in the simulation time as well as the size and distribution of the initial points to obtain an accurate but potentially coarse result. 

To demonstrate the behavior of the distances $d_H(f_1^i(\cI),f_1^{i+1}(\cI))$,  we consider a radially symmetric ODE with an attractor consisting of two periodic orbits. The domain is chosen to be the square $X = [-4,4]^2$ and  equations, in the polar coordinates, are given by
\begin{equation}\label{eqn:radial_2labels1}
    \begin{cases}
        $$\dot{\theta} = 1$$ \\
        $$\newline \dot{r} = -r(r-1)(r-2)(r-3)$$
     \end{cases}
\end{equation}
Figure~\ref{fig:hd_vary_pts}(a) shows that for large enough sample $\cI$ the distances $d_H(f_1^i(\cI) f_1^{i+1}(\cI))$ become small around $i=6$ and then they are almost constant. This suggests that the points $f_1^{I}(\cI)$ for $I> 6$ are close to the respective attractors. Figure~\ref{fig:hd_vary_pts}(b) depicts the set $f_1^{14}$  for the longest integration time of our simulation. As expected this set is a good approximation of the two periodic orbits that form the distinct attractors of the system. 

As demonstrated above, for sufficiently large $I$ and reasonably dense $\cI$, the set $f_1^I(\cI)$ is a good approximation of the attractors. To identify the approximations of the attractors $A_k$ we search for clusters in $f_1^I(\cI)$. There are many different clustering methods that can  achieve this goal~\cite{treshansky2001overview, murtagh2012algorithms}. We use a hierarchical clustering algorithm based on persistence homology~\cite{PhysD2016, persistence_cluster} of a Vietoris-Rips complex~\cite{Carlsson_2014} defined on the points in $f_1^I(\cI)$. This method is similar to clustering with dendrograms~\cite{nielsen2016hierarchical}. Points in the $0$-dimensional persistence diagram~\cite{EH08} indicate scales at which distinct clusters merge together.  The points in the persistence diagram with larger death coordinates correspond to better separated clusters~\cite{PhysD2016}.  The persistence diagram shown in Figure~\ref{fig:hd_vary_pts}(c) contains two points whose death  coordinates are considerably larger than death coordinates of the other points. This indicates that there are two distinct and well separated clusters corresponding to two circles visible in Figure~\ref{fig:hd_vary_pts}(b).  We  identify the distinct clusters with the sets $A_k$, i.e. the set $A_k$ consists of the points in $f_1^I(\cI)$ that belong to the $k$-th cluster. Finally we define the labeling function $F\colon \cI \to \{0, \ldots, K\}$ by $F(x) := k \text{ if } f_1^I(x) \in A_k$.

\subsection{Attracting neighborhood construction}
\label{sec:ML}

Potential attracting neighborhoods $N_k$ for the sets $A_k$ are constructed by using a partition of the domain $X$ into  polytopes. This partition is induced by a neural network that extends the labeling function $F\colon \cI \to \{0, \ldots, K\}$ to $X$~\cite{Balestriero2019}. There are many different machine learning methods for constructing this extension~\cite{kotsiantis2006machine}.  We consider a neural network model based on the MultiLayer Perceptron (MLP)~\cite{noriega2005multilayer, delashmit2005recent}, so that the resulting neural network model is a piecewise-linear map $F_{\theta} \colon \R^d \to \R$. Polytopes on which this map is linear depend on the network parameters (weights) $\theta$~\cite{raghu2017expressive}, and there are algorithms to compute these polytopes from  $\theta$, e.g.\ \cite{masden2022algorithmic}.

To keep the computations simple,  we use a regression perceptron network with a single hidden layer to extend the function $F$.
It is well known~\cite{Hornik, Leshno-Lin-Pinkus-Schocken} that for every continuous function $F \colon X \to \R$ defined on a compact set $X$, and $\varepsilon > 0$ there  exists a perceptron network $F_{\theta}$ with a single hidden layer and weights $\theta$ such that $\sup_{x\in X}|F(x) - F_{\theta}(x)| < \varepsilon$.  We note that our function $F$ is discontinuous at the boundaries of the basins of attraction so this result does not apply.  However, due to finite size of our sample $\cI$, there is a continuous function $\bar{F}$ that agrees with $F$ on $\cI$ and satisfies $\bar{F}(x)=k$ for every $x$ in some closed neighborhood of $A_k$ for every $k$.  

In practice, it is more efficient  to use  multiclass classification to construct the approximation of $\bar{F}$. This approach also yields  a vector function $p \colon X \to  \R^k$ whose $k$-th coordinate  $p_k(x)$ expresses the probability that $\bar{F}(x) = k$. Using this vector we could  define $N_k := \{x\in X \colon p_k(x) \approx 1\}$. If the accuracy of the classification is sufficiently high, then this set is an attracting neighborhood of $A_k$. There are many heuristic results about the accuracy of such classification  methods, but only recently have the first theoretical results  started to appear \cite{thrampoulidis2020theoretical, bianchini2014complexity}.
Moreover, this method does not in general provide a decomposition of the sets $N_k$ into polytopes.

Hence, as already mentioned, we consider a regression network with a single hidden layer of width $p$.  This network is a map $F_{\theta} \colon \R^d \to \R$ that maps the input value $x \in \R^d$ to an output prediction that depends on the weights $\theta$. In particular, $\theta$ is a triplet $(A, b$, $w)$,  where $A$ is a real matrix in $M_{p\times d}(\R)$ and $b $, $w$ are $p$-dimensional real vectors. 
The nonlinear activation function $h^{a, b} \colon \R^p \to \R^p$ is defined coordinate-wise as HardTanh
\[
h_i^{a, b}(x_i) = \left\{
	\begin{array}{ll}
        a  & \mbox{if } x_i < a, \\
		x_i  & \mbox{if } a \leq x_i \leq b, \\
		b & \mbox{if } x_i > b.
	\end{array}
\right.
\]
For a given $\theta$ the map is defined by a scalar product 
\begin{equation}
\label{eqn:net_def}
    F_{\theta}( x ) = h^{0, K-1}(\langle w, h^{0,1}(Ax + b)\rangle). 
\end{equation}

The function  $F_{\theta}$ is piecewise-linear, and we denote the supporting hyperplanes of the polytopes on which this function is linear by
$H_i^j$ with  $1 \leq i \leq p$ and $j = 0,1$. The hyperplane $H_i^j$ is defined by
\begin{equation}
\label{eqn:hyperplanes}
    \langle A[i], x \rangle + b_i = j    
\end{equation}
where $A[i]$ is the $i$-th row of the matrix $A$ and $b_i$ is the $i$-th element of the vector $b$.

To take advantage of the existing code for computing homology in higher dimensions~\cite{pyCHomP_repo}, we restrict the weights of the matrix $A$ in such a way that the relevant polytopes are parallelotopes. For simplicity, we set the width of the network  $p = qd$ for some $q \geq 1$ and divide the rows of the matrix $A$ into $d$ groups each containing $q$ rows. Within each group we require that the rows are identical.

Each group of the identical rows of the matrix $A$ induces $2q$ parallel hyperplanes defined by Equation~\eqref{eqn:hyperplanes}. Moreover, if the rows of the matrix $A$ span $\R^d$, then two hyperplanes corresponding to rows in distinct groups are not parallel. In this case, the collection of all the hyperplanes $H_i^j$ divides $\R^d$ into a collection $\{S_i\}_{i=1}^M$ of $d$-dimensional subsets that are either parallelotopes or infinitely long strips as in Figure~\ref{fig:sys1_decomp}. To ensure that the sets $S_i$ intersecting $X$ are bounded, we sometimes need to add two hyperplanes to each group of parallel hyperplanes. In the rest of the paper, the resulting collection of parallelotopes that intersect $X$ is denoted $\{C_i\}_{i=1}^N$.
Due to its cubical data structure, the decomposition of $X$ given by this method is called a Machine-Learned Cubical Decomposition (MLCD). 

We utilize the learned approximation $F_{ \theta}$ of the labeling function $F$ to define the potential attracting neighborhoods as follows.
Given $\epsilon>0,$ define
\begin{equation}\label{eqn:labeled_region}
N_k := \bigcup_{C_i\in \mathcal{N}_k}C_i \quad\text{where}\quad \mathcal{N}_k=\{C_i \colon |F_{ \theta}(x) - k | \leq \varepsilon \text{ for all } x \in C_i\}.
\end{equation}
By Bauer's minimum principle, $C_i$ is contained in  $\mathcal{N}_k$ if and only if $|F_{ \theta}(x) - k |\leq \varepsilon$ for all the vertices $x$ of $C_i$. The union of all parallelotopes $C_i$ that are not contained in any $\mathcal{N}_k$ form the \emph{uncertain region} denoted by $U$. We note that larger values of $\varepsilon$ should be used if the quality of the approximation $F_\theta$ is low.  In this paper, we do not investigate systematic methods for choosing $\epsilon$.  Instead, we consider a range of values and indicate how the results depend on $\varepsilon$ for particular problems.

We close this section by pointing out that the chosen nonlinearity $h^{a, b}$ allows us to approximate a broader class of functions than the standard nonlinearity ReLU. However, already for two-dimensional domains, this network cannot achieve an arbitrarily good approximation mentioned at the beginning of this section. To see this consider an arbitrary function $F_{ \theta},$ and let $u_1$ and $u_2$ be two linearly independent rows of the matrix $A$. Note that in the coordinates given by this basis, $F_{ \theta}(x,y) = f_1(x) + f_2(y)$ for some functions $f_1, f_2 \colon \R \to \R$. Clearly this class of functions is too restrictive to represent an arbitrary step function. Thus to study general dynamical systems it is necessary to use an unconstrained network architecture and deal with the general polytopes. In Section~\ref{sec:results} we provide heuristics for detecting if the network is sufficiently versatile to adequately approximate the labeling function. The same heuristics can be applied to unconstrained networks to decide if their width and depth is sufficient.

\subsection{Training and benchmarking}
\label{sec:training}

In this section we describe the initialization and training procedure for the neural network $F_\theta$, which is trained to approximate the labeling function $F$. The values of the training hyperparameters are reported in Section~\ref{sec:net_params}. We typically use a testing dataset with size equal to the size of the training dataset. For systems where producing balanced data is more expensive, we use a testing to training data ratio of $1:4$.

In order to ensure that data sparsity does not affect the results, we create conservatively large training datasets $\cI$. In cases where the sampling procedure produces unbalanced data, oversampling of the least witnessed class is used as the primary technique to produce balanced training data. 
We note that this procedure is solely based on the labeling produced by the clustering algorithm and does not require any knowledge of the attracting blocks.  For the ellipsoidal bistable system in four and five dimensions this method requires a fair amount of computation. We use these systems only to investigate the sizes of the grids, thus we decided to avoid the lengthy simulations of the orbits and exploit the knowledge of the attracting neighborhoods instead.
We non-uniformly sample the domain $X$ and determine to which attracting neighborhood points belong without integrating. The ratio of the points outside and inside the ellipsoid  ($70:30$ in four dimensions and $86:14$ in five dimensions) was determined from preliminary experiments. 

We did not carry out a systematic study of balancing necessary to  discover all the attractors in our systems. First of all there might be an infinite  number of attractors, and even in the finite case we suppose that required balancing will be problem dependent. So, we do not expect that  analysis of our simple examples would provide any general results or heuristics. However, even if we miss an attractor, we still produce a valid  finite lattice of attractors $\sA$ as guaranteed by  Theorem~\ref{thm:fundamental}. 

To study the performance of our method, we create $100$ different realizations of $F_\theta$ for each system.  The initial weights of each realization depend on the size of the domain $X = \prod_{i = 1}^d [a_i, b_i]$ and the number of nodes $p =qd$. We recall that we split the nodes into $d$ groups, and the $q$ nodes in each group correspond to $2q$ parallel hyperplanes in $\R^d$. We always choose the initial weights in such a way that the hyperplanes in the $i$-th group are perpendicular to the $i$-th standard basis vector $e_i$. To be more precise the rows of the matrix $A$ corresponding to the $i$-th group of nodes are set equal to the transpose of $e_i$. Only the weights corresponding to the hyperplane offsets are chosen randomly. We want these hyperplanes to be initially spread along the domain  $X$, so for the $i$-th group we sample the offsets randomly according to the uniform distribution on $(a_i, b_i)$.

To train the network, we use the Adam  optimizer~\cite{Adam} with learning rate $0.01$. We compute a test loss on a testing dataset at every epoch of training. 
We use a validation-based early stopping criterion defined by a patience parameter $\rho$ so that the  training stops if the average of the test loss over the last $\rho$ epochs stops decreasing. Otherwise we train for a problem-dependent number of epochs specified in Section~\ref{sec:net_params}.

There are two cases in which the realization is not included in the set of $100$ realizations: 1) if the normal vectors of the hyperplanes do not span $\mathbb R^d$, and 2) the training process did not converge.  To assess if the training process did not converge, we use a heuristic based on reduction of the loss between the first and last epoch.  If the fraction of the training losses in the first and last epoch is less than a given parameter $\beta$, then the training did not converge, and we do not use this realization. We do not investigate how to systematically choose $\beta$ and $\rho$, but report the values used for each example in Section~\ref{sec:net_params}.

We establish a benchmark to compare the number of polytopes from the Machine-Learned Cubical Decomposition (MLCD) to the number of grid elements needed to construct the sets $N_k$ and $U$ from a regular cubical decomposition, by which we mean a partition of $X$ by uniform subdivision in each coordinate direction, as shown in Figure~\ref{fig:elipsoidal_scaling}(c).

First, we train an (unconstrained) multilayer perceptron network using \\  {\tt sklearn.MLPClassifier} to label points according to their basin of attraction. The network is trained until the test accuracy exceeds $95\%$ on a test set consisting of $10\%$ of the original data. The vertices of a regular cubical grid on $X$ are labeled by the classifier. If all vertices of a cube have the same label $k$, then the cube is included in $N_k$; otherwise, it is included in $U$. In Section \ref{sec:results}, we consider the minimum grid size required to recover correct Morse representations and Conley indices of the corresponding Morse sets, as described by the homology of the sets $N_k$ and $U$ listed in Section~\ref{sec:indexes}.

\section{Attractor Identification Results}
\label{sec:results}

In this section we demonstrate that the  Machine-Learned Cubical Decomposition (MLCD) can provide correct Morse representations of a dynamical system together with Conley indices of its Morse sets. We also show that the number of polytopes in the MLCD is significantly smaller than the number of cubes required by a regular cubical decomposition. In Section~\ref{sec:fixed_points} we start with simple dynamical systems whose attractors consist of fixed points and familiarize the reader with the workflow of MLCD. Then in Section~\ref{sec:periodic_orbits} we move to systems containing periodic orbits and highlight the small number of polytopes in MLCD. Finally, in  Section~\ref{sec:complex} we consider more complex systems.

All the statistics presented are based on $100$ independent realizations of the network $F_\theta$, described in Section~\ref{sec:training}. While we cannot achieve the rigor of the universal approximation theorem~\cite{Hornik}, we  argue that the following remark, corroborated by all the examples considered in this paper, provides a useful guideline for practical applications.

\begin{remark}
\label{rem:heuristics}
If the architecture of the network is sufficiently expressive and a sufficiently small value of the loss function is achieved during the training, then the approximation $F_\theta$ of the labeling function $F$ should lead to the correct Conley indices. Moreover, if the computed Conley indices are the same for all the realizations with the final testing loss below some value, then this value is likely to be sufficiently small. However,  if the  indices keep changing for all observed loss values, then the  network is not sufficiently expressive. 
\end{remark}
The rational of this remark relies on the fact that if the  testing loss is small, then $F_\theta$ closely approximates the piecewise constant labeling function $F$. Hence it is likely that for small values of $\varepsilon$, each set $N_k$, defined by~\eqref{eqn:labeled_region}, is a subset of $F^{-1}(k)$ and has the same Conley index as the attractor $A_k$. In this paper, we consider $\epsilon \in \{0.1, 0.2, 0.3, 0.4, 0.49\}$ and often report the results for the  $\varepsilon$ that yields the highest success rate. The detailed study of the dependence on $\varepsilon$ is beyond the scope of this paper, and we only discuss it briefly in Section~\ref{sec:conclusion}.

\begin{figure}[]
    \centering
    \includegraphics[width=13cm]{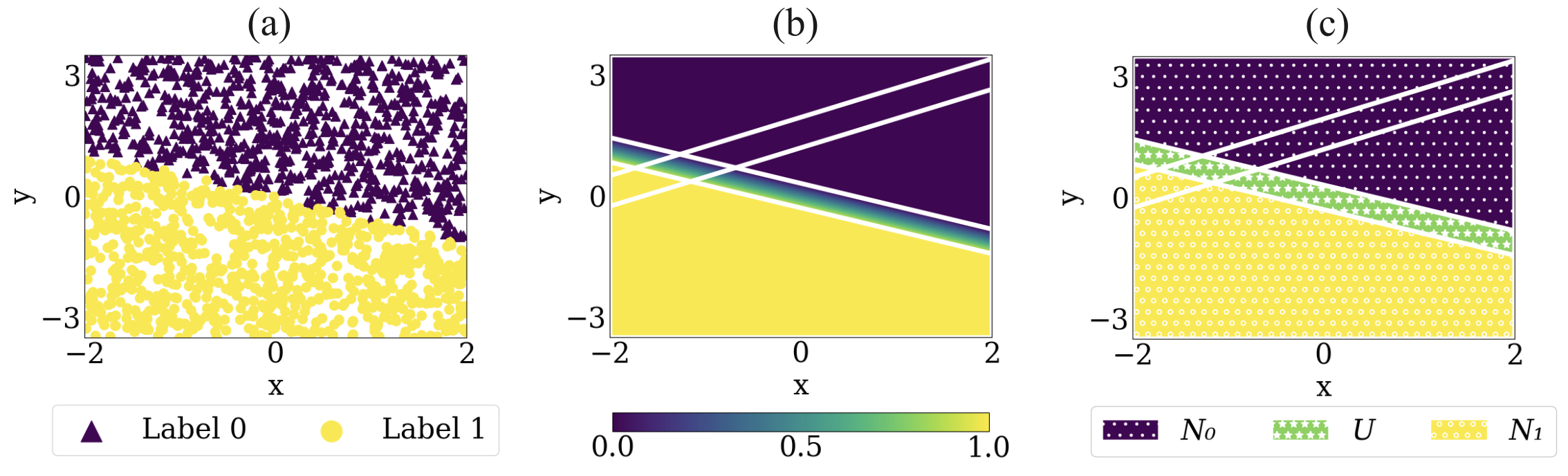}
    \caption{(a)~A sample of the labeling function $F$ for the bistable, two-dimensional system given by~\eqref{eqn:sys1}. (b)~$F_\theta$ (single realization) approximates $F$ and is linear on the parallelotopes bounded by the hyperplanes displayed in white. These hyperplanes partition the set $X$ into one parallelotope and eight infinitely long strips. By adding two hyperplanes to each set of the parallel hyperplanes,  shown in white, the set $X$ is partitioned into nine finite parallelotopes.  (c)~Parallelotopes in the sets $N_0$, $N_1$ and $U$ for the given $F_\theta$. In this case the same parallelotopes are obtained for all considered values of $\varepsilon$.  }
    \label{fig:sys1_decomp}
\end{figure}

\subsection{Fixed points}
\label{sec:fixed_points}
The first system that we consider is an ODE on the rectangle $X=[-2,2] \times [-3.5,3.5]$ 
given by:
\begin{equation}\label{eqn:sys1} 
\dot{x} = u(1-u^2)\quad\text{and}\quad \dot{y} = u^2 (3-2u^2) - v     
\end{equation}
where
\begin{equation*}
 u=\frac{x}{2} + \frac{\sqrt{3}}{2}y  \quad\text{and}\quad v=\frac{y}{2} -\frac{\sqrt{3}}{2}x. 
\end{equation*}
This system has two stable fixed points at approximately $(-0.366,1.366)$ and $(-1.366,-0.366)$. Their basins of attraction are separated by the straight line $y=-\sqrt{3}x$, which is the stable manifold of the saddle point at the origin.

To produce the MLCD we first sample the system and construct the labeling function $F \colon \cI \to \{0,1\}$ as described in Section~\ref{sec:AI}. The shape of the function $F$, shown in Figure~\ref{fig:sys1_decomp}(a), clearly indicates two distinct basins of attractions separated by a straight line. Due to this simple geometry a very small network provides a good approximation of $F$. Figure~\ref{fig:sys1_decomp}(b) documents that one hidden neuron (two hyperplanes) per dimension is sufficient to separate the basins of attraction from the separatrix. Hence, we use a network with width $p=2$. 

For two-dimensional input, the constrained network $F_\theta$ with two hidden neurons  typically  divides $\mathbb R^2$ into nine subsets out of which eight are unbounded.  For the $F_\theta$ depicted in Figure~\ref{fig:sys1_decomp}(b), we need to add two hyperplanes in each direction to ensure that the domain $X$ is covered by finite parallelotopes. After adding the hyperplanes, the set $X$ is covered by nine  parallelotopes and the sets $N_0$, $N_1$ and $U$, shown in Figure~\ref{fig:sys1_decomp}(c), are identical for all considered values of $\varepsilon$. Moreover, their Conley indices are correct; see Section~\ref{sec:indexes} for the correct Conley indices of all the systems considered in this paper.

\begin{figure}
    \centering
    \includegraphics[width=\linewidth]{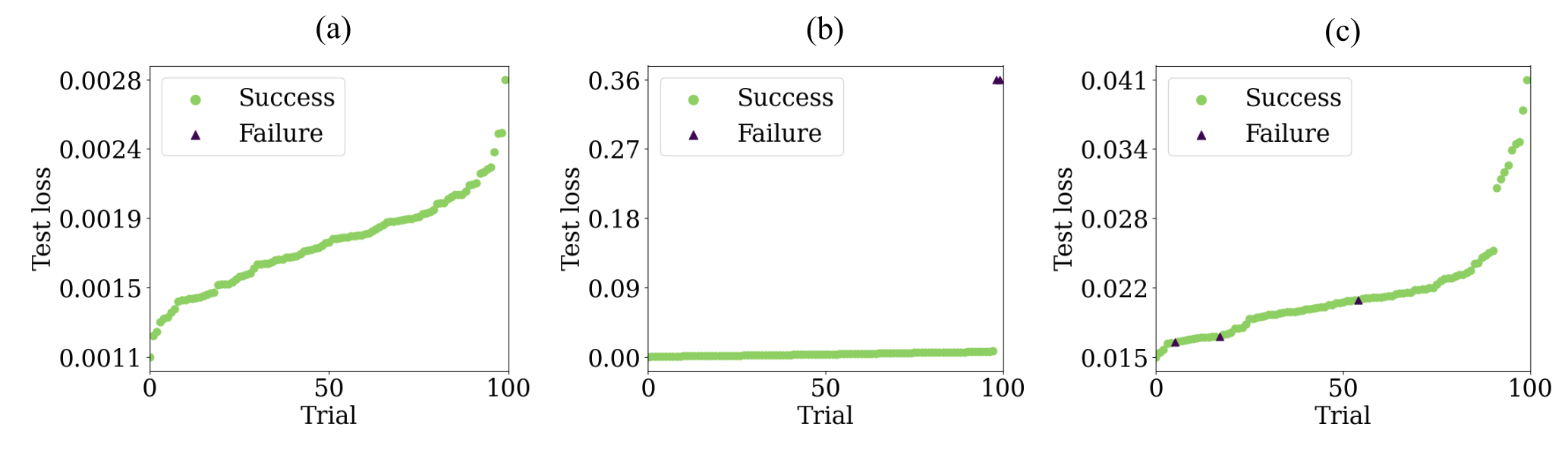}
    \caption{The final values of the test loss  for successes (green points) and failures (violet triangles) obtained by networks with one node per dimension for the (a) Two-dimensional system with the linear separatrix given by~\eqref{eqn:sys1}, (b) Four-dimensional system with the nonlinear separatrix given by~\eqref{eqn:curved_separatrix}, and (c) Six-dimensional EMT Hill system given by~\eqref{eqn:EMT}.
    }
    \label{fig:sys_with_fixed_points}
\end{figure}

Now we study the performance on $100$ independent realizations of $F_\theta$. For each realization we record the final test loss and the number of polytopes in the MLCD (by definition, all polytopes intersect $X$). We also check if the Conley indices computed using the MLCD are correct. If yes, we label the realization as a success. Otherwise we label it as a failure. Figure~\ref{fig:sys_with_fixed_points}(a) shows the final test losses for the  individual realizations. All these values are relatively small and the sets  $N_0, N_1$ and $U$ constructed from grids generated by each realization $F_\theta$ always yield the correct Conley indices. The average number of polytopes (over $100$ trials) is six while the standard deviation is one. Hence, even in this simple case MLCDs slightly outperform the regular cubical decomposition, which requires nine uniform cubes to obtain the correct Conley index. Moreover, the number of polytopes is insensitive to the  orientation of the separatrix while the number of cubes increases significantly if the separatrix  is closely aligned with the diagonal of the domain.

Next we consider a higher-dimensional system with two stable fixed points in the domain $X = [-2, 2] \times [-3.5, 3.5] \times [-2, 2] \times [-2, 2]$. The system is given by the following ODE:
    \begin{equation}\label{eqn:curved_separatrix}
         \begin{cases}
             $$\dot{x_1} = -x_1$$ \\
             $$\dot{x_2} = (x_2-x_1^2)(9-x_2^2)$$\\
             $$\dot{x_3} = -x_3$$\\ 
             $$\dot{x_4} = -x_4$$
         \end{cases}
     \end{equation}
The basins of attraction of the two stable fixed points are now separated by a three-dimensional stable manifold of a saddle equilibrium. 
Since neither the curvature of the manifold nor the size of the set $X$ is too large, it is again possible to delineate the separatrix from the basins of attraction by using a neural network with one node per dimension.

In this case, $2\%$ of the trained networks are failures while the other $98\%$ are successes for all values considered of $\varepsilon$.  Figure~\ref{fig:sys_with_fixed_points}(b) shows the final test loss for all realizations from both successes and failures. All the trials with small final loss lead to the same Conley indices. As predicted by  Remark~\ref{rem:heuristics}, all these trials are successes. The improvement over the regular cubical decomposition, which requires $81$ cubes, is more remarkable because on average MLCD only produces $33$ parallelotopes (the standard deviation is $10$).  

We also test the very small network architecture ($q=1$ and $p=6$) on a six-dimensional system defined by~\eqref{eqn:EMT}. This system contains two fixed points, but the geometry of the separatrix is more complicated. We still obtain a $97\%$ success rate. The small size of the network 
results in a mean number of 331 parallelotopes in the MLCD (with a standard deviation of 90), whereas the regular cubical decomposition produces 46,656 cubes. However, the final test losses, depicted in Figure~\ref{fig:sys_with_fixed_points}(c), do not indicate stabilization of the Conley indices for small values.  As predicted by Remark~\ref{rem:heuristics}, the network is too small to sufficiently approximate the labeling function. This is also corroborated by the relatively high value of the lowest loss achieved over $100$ realizations. We note that in this case, the results are very sensitive to the particular choice of $\varepsilon$. When the size of the network is increased to
$q=2$ nodes per dimension and width $p=12$, the separation between the successes and failures is more pronounced and the results are less sensitive to the value of $\varepsilon$. Even for this larger network, the average number of parallelotopes is $2,631$ which is still a significant improvement. However, we do not achieve a convincingly low test loss. This is likely a consequence of the restriction imposed on the matrix $A$ in Section~\ref{sec:ML}. We discuss this problem in more detail in Section~\ref{sec:complex}.

\subsection{Fixed points and periodic orbits} 
\label{sec:periodic_orbits}

\begin{figure}
    \centering
    \begin{subfigure}[c]{0.4\textwidth}
        \centering
        \includegraphics[width=5cm]{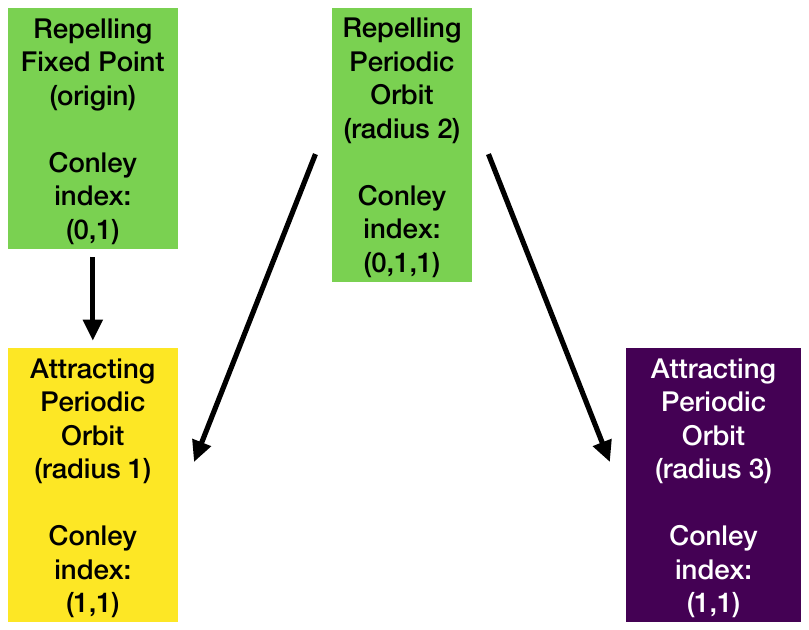}
        \caption{ }
    \end{subfigure}
    \qquad\qquad
    \begin{subfigure}[c]{0.4\textwidth}
        \centering
        \includegraphics[width=5cm]{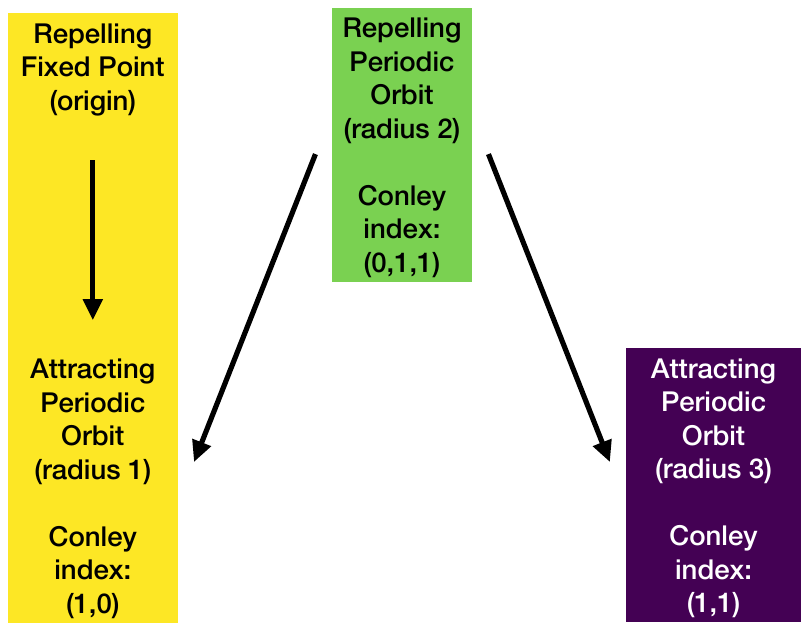}
        \caption{ }
    \end{subfigure}
    \caption{(a) Finest Morse representation for the radially symmetric system defined by Equation~\eqref{eqn:radial_2labels}. (b) A coarsened Morse representation that still exhibits the bistability of the system.}\label{fig:MR_radial}
\end{figure}

\begin{figure}
    \centering
    \includegraphics[width=0.95\linewidth]{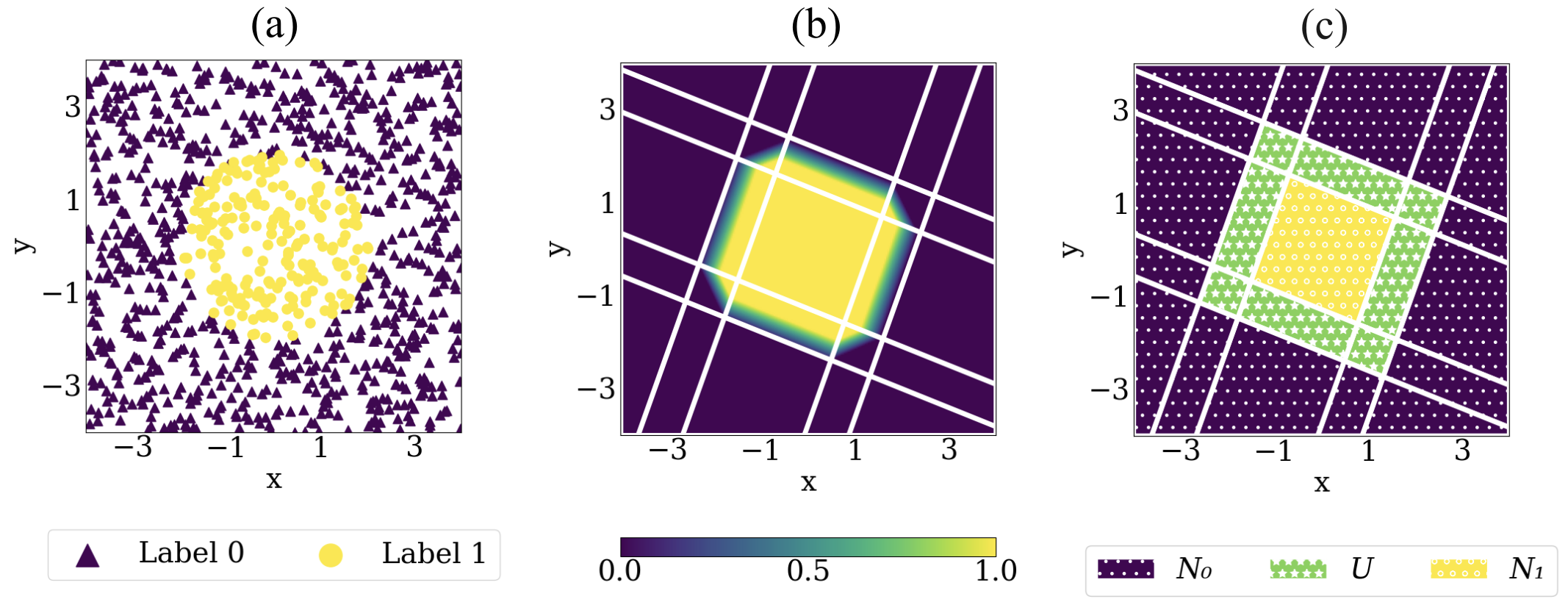}
    \caption{(a) A sample of the labeling function $F$ for the radially symmetric system defined by Equation~\eqref{eqn:radial_2labels}. (b) $F_\theta$ (single realization) approximates $F$ and is linear on the parallelotopes bounded by the hyperplanes displayed in white. (c) Parallelotopes in the sets $N_0$, $N_1$ and $U$ for the given $F_\theta$.  }
    \label{fig:2d_radial}
\end{figure}

In this section we consider systems with basins of attraction of more complicated structure.
The first system is an ODE on the set $X=[-4,4] \times [-4,4]$. The equations are  given in polar coordinates by:
    \begin{equation}\label{eqn:radial_2labels}
        \begin{cases}
            $$\dot{\theta} = 1$$ \\
            $$\newline \dot{r} = -r(r-1)(r-2)(r-3)$$
        \end{cases}
    \end{equation}
The asymptotic dynamics of the system includes two attracting periodic orbits, a repelling periodic orbit as a separatrix, and the origin as a repelling fixed point. Figure~\ref{fig:MR_radial}(a) shows the finest Morse representation. In this representation, the repelling fixed point at the origin is a distinct Morse set. However, due to the finite nature of the data, this unstable structure cannot be fully resolved. The constructed labeling function depicted in  Figure~\ref{fig:2d_radial}(a) misses the fact that the origin is a fixed point. Hence, the neural network $F_\theta$ presented in Figure~\ref{fig:2d_radial}(b) combines the origin, the attracting limit cycle of radius one, and the connecting orbits from the origin to the periodic orbit into a single attractor. Figure~\ref{fig:MR_radial}(b) shows a coarser Morse representation obtained for the sets $N_0, N_1,$ and $U$ depicted in Figure~\ref{fig:2d_radial}(b). Although we cannot find the finest Morse representation,  we still properly capture the bistability in the system.

Figure~\ref{fig:2d_radial} demonstrates that two hyperplanes per dimension are not sufficient to separate the basins of attraction and the separatrix. We use the next smallest network and set the number of hidden nodes per dimension equal to $q=2$, which results in a network of width $p=4$. The same figure shows that this small network can provide a sufficient approximation $F_\theta$ of the labeling function. In particular, the sets $N_1, N_2,$ and $U$ yield the correct Conley indices of the Morse sets in the coarser representation.

\begin{figure}
    \centering
    \includegraphics[width=\linewidth]{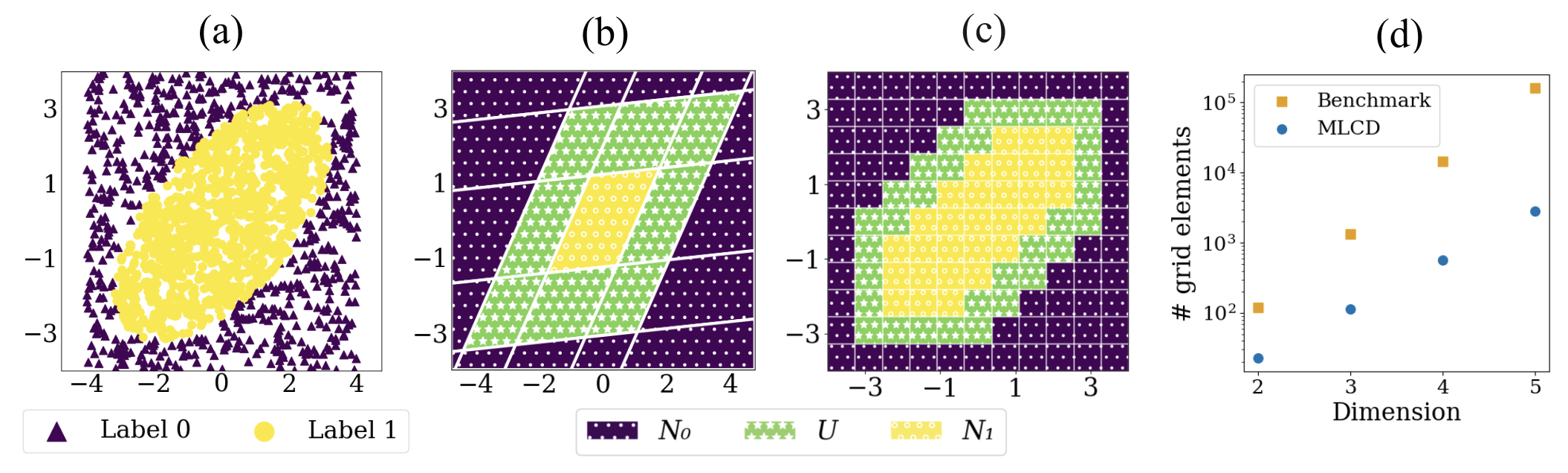}
    \caption{(a)~A sample of the labeling function $F$ for the radially symmetric system defined by Equation~\eqref{eqn:radial_2labels} after the coordinate change. (b)~Parallelotopes produced by the MLCD that form the sets $N_0$, $N_1$ and $U$. (c)~Cubes in the regular cubical decomposition that form the sets $N_0$, $N_1$ and $U$. d)~Number of parallelotopes and cubes necessary to capture the sets $N_0$, $N_1$ and $U$ versus the dimension of the system.}
    \label{fig:elipsoidal_scaling}
\end{figure}

More complicated geometry of the attractor and the small size of the network result in a lower success rate, which is $31\%$. All the realizations with sufficiently small final test loss yield the same indices and as predicted by Remark~\ref{rem:heuristics}, they are correct.  However, for this system there is no improvement over the uniform grid. Figures~\ref{fig:2d_radial}(b)~and~(c) show a typical grid obtained by  MLCD and the structure of the sets $N_0, N_1$ and $U$. There are always $25$ polytopes which is exactly the necessary number of cubes in the regular grid.
Figure~\ref{fig:elipsoidal_scaling} shows that if the system is transformed so that the periodic orbit is an ellipse with diagonal-aligned axes, then the improvement is significant.

More precisely, to study the dependence of the number of grid elements on the dimension we modify the system defined by Equation~\eqref{eqn:radial_2labels} by a coordinate change. This coordinate change deforms the circular periodic orbit to an ellipsoidal periodic orbit. The labeling function (see Figure~\ref{fig:elipsoidal_scaling}(a)) shows the geometry of the basins of attraction for the modified system. We extend this system to $3,4$ and $5$ dimensions so that the separatrix given by the unstable periodic orbit changes to an ellipsoid that separates the space into two basins of attraction. 

\begin{figure}
    \centering
    \includegraphics[width=\linewidth]{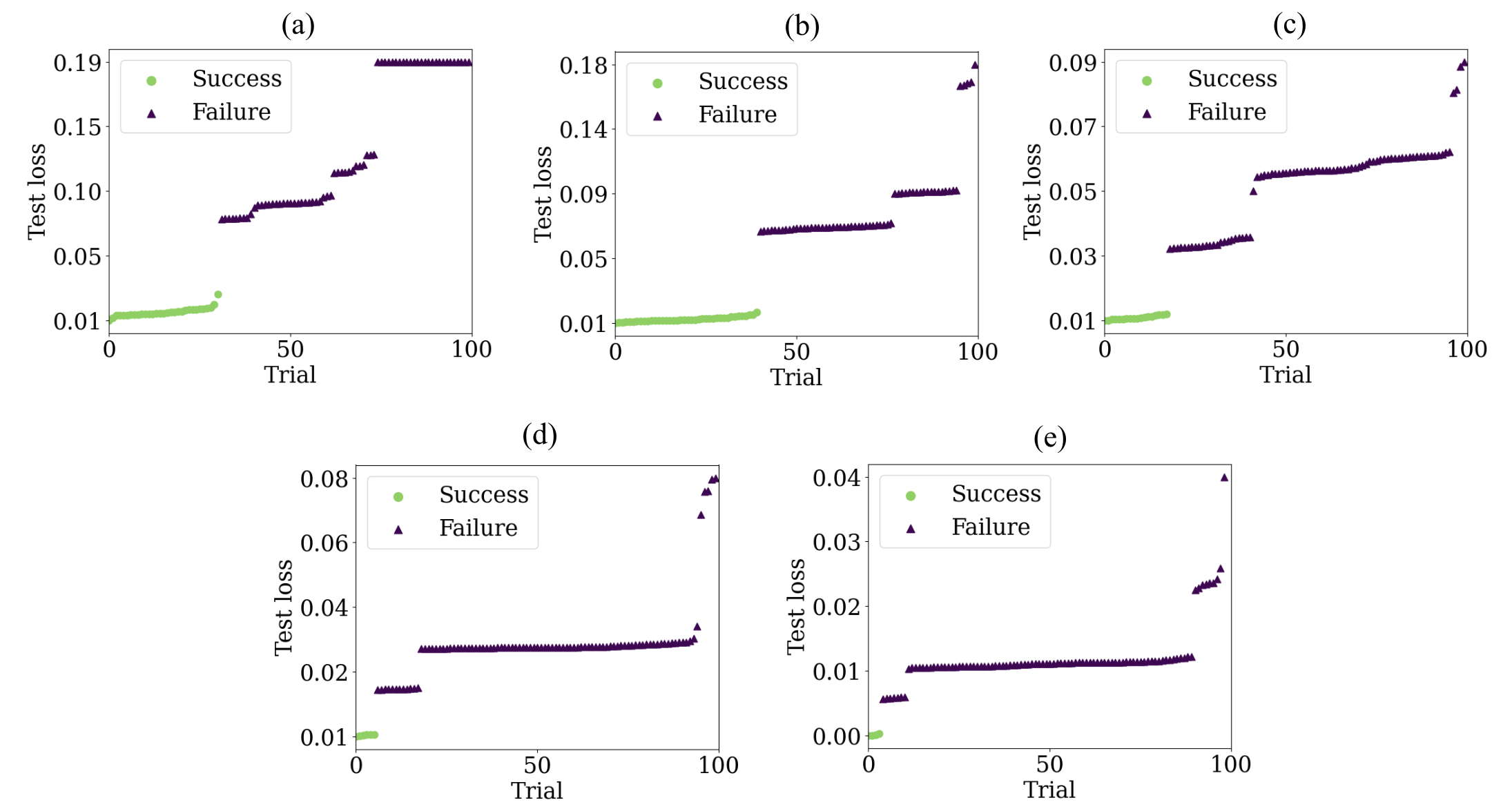}
    \caption{The final values of the test loss  for successes (green points) and failures (violet  triangles)  obtained by networks with two nodes per dimension.  (a) Two-dimensional radial bistable system given by~\eqref{eqn:radial_2labels}. The ellipsoidal bistable system in dimension (b) 2, (c) 3, (d) 4 and (e) 5.}
    \label{fig:dist_radial}
\end{figure}

To approximate the labeling function for the higher-dimensional systems we again use the network architecture with two hidden nodes per dimension. Figure~\ref{fig:dist_radial} indicates that as predicted by Remark~\ref{rem:heuristics},  the computed Conley indices are correct for a sufficiently small test loss. We note that realizations that achieve a reasonably small test loss become progressively more rare as the dimension of the systems increases. This leads to a progressively smaller percentage of successes. In particular, the success rates are $40\%, 18\%, 6\%$ and $4\%$ as the dimension increases. However, training the network is inexpensive, and a realization with a small loss can be found very quickly. In what follows we demonstrate that the amount of time required to find a low loss approximation is negligible compared to the savings in computing the Conley indices (homology).

Figure~\ref{fig:elipsoidal_scaling}(b) shows a typical grid produced by MLCD. In each direction, the set $X$ intersects with $4$ hyperplanes. Hence, the MLCD of $X$ contains $5^2$ parallelotopes. Figure~\ref{fig:elipsoidal_scaling}(c) depicts the regular grid, formed by subdividing the set $X$ with $10$ hyperplanes in each direction, resulting in $11^2$ cubes.  As we increase the dimension, the number of required hyperplanes per dimension is roughly preserved, i.e. $4$ for MLCD and $10$ for the uniform grid. Therefore, the number of grid elements grows exponentially for both grids, but the base of the exponential function is reduced by more than half by using MLCD. Figure~\ref{fig:elipsoidal_scaling}(d) indicates that as the dimension $n$ of the system increases, the grid sizes indeed grow as $5^n$ and $11^n$ respectively.

\subsection{Systems with more complex attractors} 
\label{sec:complex}

\begin{figure}
    \centering
    \includegraphics[width=\linewidth]{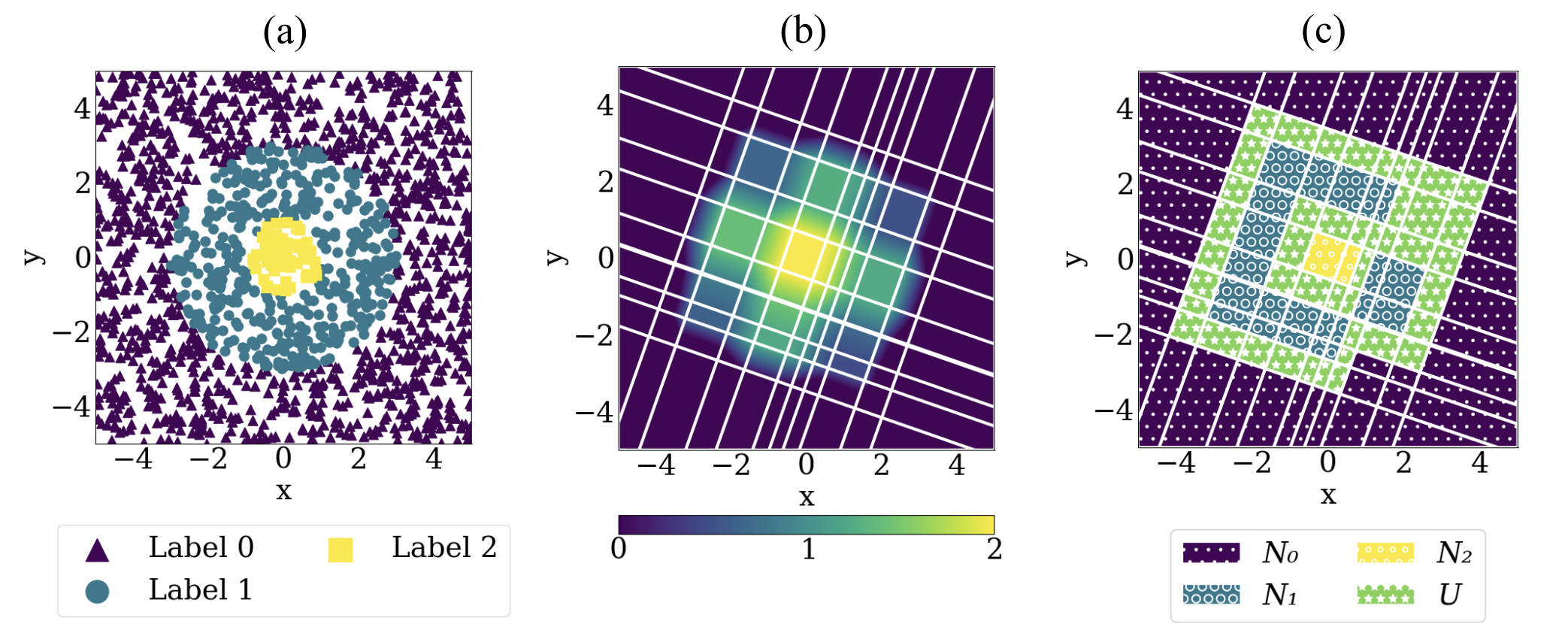}
    \caption{(a)~A sample of the labeling function $F$ for the system defined by Equation~\eqref{eqn:radial_3labels}. (b)~$F_\theta$ (single realization) demonstrates that the restricted network architecture cannot properly approximate the labeling function.  Thus (c)~the sets $N_0$, $N_1, N_2,$ and $U$ do not capture the radially symmetric nature of the basins of attraction and separatrixes.}
    \label{figures:radial_3failure}
\end{figure}

In Section~\ref{sec:ML}, we describe a restriction of the weight matrix $A$ of the neural network \eqref{eqn:net_def} so that the polytopes in the grid are parallelotopes. This restriction allows us to easily identify the grid elements and also to use efficient (cube-based) code for computing homology.  In this section, we present two simple examples to show that this approach is too restrictive.  The first system is given by an ODE
    \begin{equation}\label{eqn:radial_3labels}
        \begin{cases}
            $$\dot{\theta} = 1$$ \\
            $$\newline \dot{r} = -r(r-1)(r-2)(r-3)(r-4)$$
        \end{cases}       
    \end{equation} 
on the domain $X=[-5,5] \times [-5,5]$. This system has two attracting periodic orbits, an attracting fixed point at the origin, and two repelling periodic orbits acting as separatrixes.  Figure~\ref{figures:radial_3failure}(a) depicts the sample of the labeling function for this system. 

\begin{figure}
    \centering
    \includegraphics[width=0.8\linewidth]{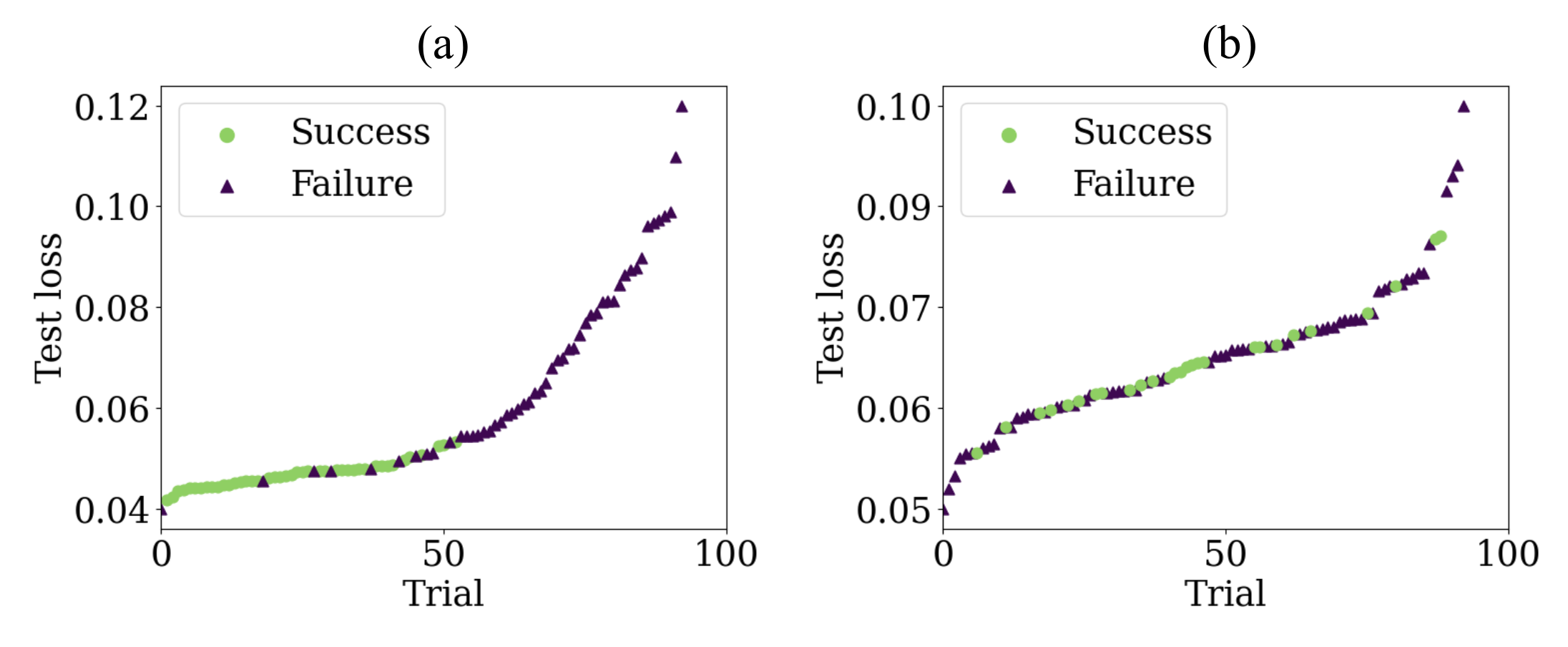}
    \caption{The final values of the test loss  for successes (green points) and failures (violet triangles) for  (a) two-dimensional system given by~\eqref{eqn:radial_3labels} and (b) three-dimensional Hill system defined by~\eqref{eqn:periodic_3d}.}
    \label{fig:dist_radial_failure}
\end{figure}

We recall that the restricted network can be expressed as $F_\theta(x,y) = f_1(x) + f_2(x)$ by using a coordinate system given by the basis of the independent rows of the matrix $A$, which is aligned with the hyperplanes (white lines) shown in Figure~\ref{figures:radial_3failure}(b).
Regardless of the size of the network, it is impossible to find a function $F_\theta$ that properly approximates the radially symmetric labeling function with three distinct integer values. Typically, the region where the labeling function is equal to one is poorly resolved, and subsequently the set $N_1$ often does not have the right topology (homology).  Figure~\ref{fig:dist_radial_failure}(a) shows that the correct homology is sometimes recovered. However, the large value of the lowest test loss and the fact that there is no separation between the successes and failures suggest that the architecture of the network is not sufficiently expressive and the results should not be trusted.

For the three-dimensional Hill system defined by  Equation~\eqref{eqn:periodic_3d}, the labeling function $F$ and its approximation cannot be easily visualized. However, we know that this system has four minimal attractors consisting of a periodic orbit and three fixed points.  Even without plotting the basins of attraction, we can detect that our method is not successful.  This  again follows from Remark~\ref{rem:heuristics}. Inspection of the final loss values, plotted in Figure~\ref{fig:dist_radial_failure}(b), reveals that they are all relatively large and the Conley indices do not stabilize as we approach the minimal loss achieved over the $100$ runs.

This example highlights one more issue that has to be carefully considered. In particular,  our simulation suggests that one fixed point has a comparatively small basin of attraction; points limiting to this fixed point make up less than $4\%$ of the initial dataset $\cI$. If the data is not properly balanced, then training the network can result in small loss even if one of the minimal attractors is missed.  The last two examples clearly demonstrate that further research and algorithm development is necessary to realize the full potential of the proposed method.

\section{Conclusion}\label{sec:conclusion}

\begin{table}
\small 
    \centering
    \begin{tabular}{|p{3.11cm}|r|r|r|r|r|}
 \hline
    \multicolumn{6}{|c|}{Comparison of Cube Counts to Benchmark} \\
 \hline
    \multicolumn{3}{|c|}{} & \multicolumn{1}{|c|}{Regular grid} & \multicolumn{2}{|c|}{MLCD} \\
 \hline
 Example & Eqn.~\# & dim($X$) & Min. \# Cubes & Mean \# cubes & SDEV\\
 \hline

 Linear separatrix & \eqref{eqn:sys1} & 2 & 9 & 6 & 1\\
 \hline
 Radial bistable & \eqref{eqn:radial_2labels} & 2 & 25 & 25 & 0\\
  \hline
 Radial tristable & \eqref{eqn:radial_3labels} & 2 & 169 & 120 & 3\\
  \hline
 Nonlinear separatrix & \eqref{eqn:curved_separatrix} & 4 & 81 & 33 & 10 \\
  \hline
   Hill system with PO & \eqref{eqn:periodic_3d} & 3 & 1,728 & 705 & 177 \\
  \hline
 EMT Hill system & \eqref{eqn:EMT} & 6 &  46,656 & 331 & 90\\
  \hline 
 Ellipsoidal bistable & Sec.~\ref{sec:periodic_orbits} & 2 & 121 & 23 & 0\\
  \hline
 Ellipsoidal bistable & Sec.~\ref{sec:periodic_orbits}& 3 & 1,331 & 115 & 2\\
  \hline
 Ellipsoidal bistable & Sec.~\ref{sec:periodic_orbits} & 4 & 14,641 & 571 & 1\\
  \hline 
 Ellipsoidal bistable & Sec.~\ref{sec:periodic_orbits} & 5 & 161,051 & 2,852 & 45\\
 \hline
\end{tabular}

\caption{Comparison of the minimum number of cubes needed to obtain the correct Conley indices using a regular cubical decomposition to the mean and standard deviation (SDEV) of the number of cubes from successful trials of MLCD. }
\label{tab:comparison_table}
\end{table}

In this paper we presented a data-driven method for quantifying the dynamics of complex systems. The method is based on Conley's topological approach that provides a framework for quantifying the global dynamics. The main advantage of this framework over a wide range of other methods is its robustness with respect to perturbations and hence to noise. It has been successfully applied to a variety of systems. However, its application to high-dimensional systems is still challenging. The major computational bottleneck is the discretization of the phase space, which is typically done using a regular cubical decomposition. The number of grid elements in this grid grows exponentially with the dimension. While the exponential growth of the grid elements cannot be avoided, we demonstrated that it is possible to dramatically reduce their number by using machine-learned cubical decompositions. Table~\ref{tab:comparison_table} highlights the reduction in the number of the grid elements necessary to properly resolve the dynamics for the systems considered in Section~\ref{sec:results}. 
 
Most data-driven methods do not provide any theoretical guarantees that the reconstructed dynamics is correct.  For our method, the main limitation for providing rigorous theoretical guarantees stems from  a lack of control over the quality of the approximation $F_\theta$ of the labeling function $F$.  Assessing the quality of the approximation $F_\theta$ produced by a neural network is a deep open problem in machine learning.
In most applications, quality of the approximation is estimated by holdout methods in which the loss function is evaluated on a part of the dataset that is not used for training.
We used this final testing loss in Remark~\ref{rem:heuristics} to formulate the following heuristic. We claim that if the final testing loss is sufficiently small then the computed Conley indices are correct. This is corroborated by all the examples considered in Section~\ref{sec:results}. To decide if the final loss is sufficiently small we compute several realizations of the function $F_\theta$. If the computed Conley indices are the same for all the realizations with the final testing loss below some value, then this value is likely to be sufficiently small.

To carry out the above mentioned computations, we have restricted the weights of the neural network \eqref{eqn:net_def}. As discussed in Section~\ref{sec:ML}, this restriction ensures that a grid decomposition consists of parallelotopes, and we can use efficient algorithms based on a cubical complex to compute homology in higher dimensions \cite{pyCHomP_repo}.  However, as discussed in Section~\ref{sec:complex}, this approach is too restrictive, and techniques based on more general polygonal decompositions are necessary to efficiently handle more complex dynamics.  In \cite{masden2022algorithmic}, methods to extract the polygonal decomposition generated by a ReLU neural network are presented for dimensions $2$ and $3$, but the efficacy and efficiency of these methods in higher dimensions require further investigation.

Algorithms for computing the homology of a general polygonal complex are presented in \cite{harker:mischaikow:mrozek:nanda} and implemented in \cite{pyCHomP_repo}.
Nevertheless, the implementation in the context of cubical complex benefits from the rigid cell structure of cubes \cite{harker2021morsetheoretictemplateshigh} that allows for efficient computations even in moderately high dimensions, e.g., ten dimensions \cite{harker:mischaikow:spendlove}. To successfully deal with higher dimensional problems requires new ideas. A potentially promising approach (see \cite{vieira123}) is to use autoencoding neural networks to encode the dynamics in a lower-dimensional latent space on which the techniques of this paper can be applied.

Understanding global dynamics requires identifying unstable structures as well as attractors. However,  
for all examples in Section~\ref{sec:results}, we do not attempt to resolve the structure of the global dynamics beyond the minimal attractors, i.e.\ we do not resolve the unstable Morse set $M$ shown in Figure~\ref{fig:attMorse}(c). Indeed, we group all cells of the MLCD that do not belong to an approximate attracting neighborhood $N_k$ into a single uncertain set $U$. 

In Section~\ref{sec:complex}, the dynamics of the system given by~\eqref{eqn:radial_3labels} has three attractors separated by two unstable periodic orbits. As shown in Figure~\ref{fig:dist_radial_failure}, the MLCD often does recover the correct homology, but there is no separation between successes and failures. However, in the successful cases, the uncertain set $U$ has two connected components (and has the homology of two disjoint circles). 

\begin{figure}
    \centering
    \includegraphics[width=0.45\linewidth]{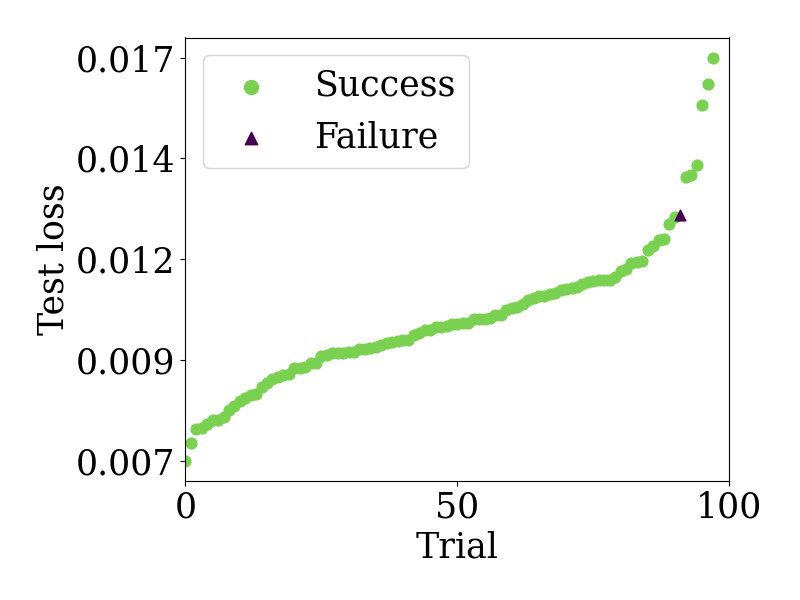}
    \caption{The final values of the test loss  for successes (green points) and failures (violet triangle) for a three-dimensional system with two repelling fixed points and three attracting fixed points.}
    \label{fig:M-loss}
\end{figure}

Figure~\ref{fig:M-loss} shows the final values of the test loss for a three-dimensional system with three attracting fixed points and two unstable fixed points. 
In this case, success is characterized by whether $U$ has two connected components,
and the existence of three attracting neighborhoods
each with Conley index $(1,0,0,0)$. There is success for all trials with sufficiently low loss, and there is separation between the successes and failures, since the geometry of the sets is more simple. This suggests that these methods can be used to characterize further unstable structures, and we leave this to future research. 

\section*{Acknowledgements}
The authors acknowledge the Office of Advanced Research Computing (OARC) at Rutgers, The State University of New Jersey for providing access to the Amarel cluster and associated research computing resources that have contributed to the results reported here. GitHub Copilot (2024) was used with duplicate detection filtering set to ``block" while developing some of the code that supplements the paper. 

\printbibliography

\newpage

\clearpage
\markboth{APPENDIX}{APPENDIX}
\section{Appendix}\label{sec:appendix}

\subsection{Parameters and Hyperparameters for neural networks}
\label{sec:net_params}

In Section \ref{sec:training}, we describe how the networks $F_\theta$ are trained. For each considered example, Table \ref{tab:example_params} lists the network width ($p$), the number of training points ($\# \cI$), the batchsize, number of epochs, patience parameter ($\rho$) for early stopping, and loss reduction threshold ($\beta$) for determining convergence.

\begin{table}[h]
    \smaller
    \centering
    \begin{tabular}{|l|r|r|r|r|r|r|r|}
    \hline
    \multicolumn{8}{|c|}{Table of Training Parameters and Hyperparameters} \\
    \hline
    \hline
    Example & dim($X$) & $p$ & $\# \cI$ & batchsize & $\#$ epochs & $\rho$ & $\beta$ \\
    \hline
    Linear separatrix & $2$ & $2$ & $10^4$ & $10^3$ & $10^2$ & $10$ & $10^{-1}$ \\
    \hline
    Radial bistable & $2$ & $4$ & $10^3$ & $10^2$ & $10^3$ & $10^2$ & $10^{-1}$  \\
    \hline
    Radial tristable & $2$ & $10$ & $10^3$ & $10^2$ & $10^3$ & $10^2$ & $5 \cdot 10^{-1}$ \\
    \hline
    Nonlinear separatrix & $4$ & $4$ & $10^4$ & $10^3$ & $10^2$ & $20$ & $10^{-1}$ \\
    \hline
     Hill system with PO & $3$ & $15$ & $1.1 \cdot 10^4$ & $10^3$ & $10^3$ & $10^2$ & $5 \cdot 10^{-1}$ \\
    \hline
    EMT Hill system& $6$ & $12$ & $10^5$ & $10^4$ & $2 \cdot 10^3$ & $10^2$ & $10^{-1}$\\
    \hline
   
    Ellipsoidal bistable & 2 & $4$ & $10^4$ & $10^3$ & $10^3$ & $10^2$ & $10^{-1}$\\
    \hline
    Ellipsoidal bistable & 3 & $6$ & $2 \cdot 10^4$& $2 \cdot 10^3$ & $10^3$ & $10^2$ & $10^{-1}$ \\
    \hline
    Ellipsoidal bistable & 4 & $8$ & $8 \cdot 10^5$ & $8 \cdot 10^4$ & $2 \cdot 10^3$ & $10^2$& $10^{-1}$\\
    \hline
    Ellipsoidal bistable & 5 & $10$ & $8 \cdot 10^5$ & $8 \cdot 10^4$ & $2 \cdot 10^3$ & $10^2$ & $10^{-1}$\\
    \hline
    \end{tabular}
    \caption{Parameters and hyperparameters used in training realizations of $F_\theta$ as described in Section \ref{sec:training}.}\label{tab:example_params}
\end{table}

\subsection{Conley indices of the attractors for the considered examples}
\label{sec:indexes}
For the examples in this paper we know a priori the Conley indices of the attractors $A_k$.
These are presented in Table~\ref{tab:correct-homology}.
We compute the homology of each $N_k$ (see \eqref{eqn:labeled_region}) to determine if we have identified the correct Conley indices.
We similarly check for the correctness of  the homology for the uncertain region $U$.

\begin{table}[h]
    \smaller
    \centering
    \setlength{\tabcolsep}{3pt} 
    \renewcommand{\arraystretch}{0.9}
    \begin{tabular}{|l|c|c|c|c|c|c|}
     \hline
     \multicolumn{7}{|c|}{Table of Correct Betti Numbers $(\beta_0,\beta_1,\ldots, \beta_d)$} \\
     \hline
     \hline
     Example & $\text{dim}(X)$ & $U$ & $N_0$ & $N_1$ & $N_2$ & $N_3$\\
     \hline
     Linear separatrix & 2 & $(1, 0, 0)$ & $(1, 0, 0)$ &  $(1, 0, 0)$ & & \\
     \hline
     Radial bistable & 2 & $(1, 1, 0)$ & $(1, 1, 0)$ & $(1, 0, 0)$ & & \\
      \hline
     Radial tristable & 2 & $(2, 2, 0)$ & $(1, 1, 0)$ & $(1, 1, 0)$ & $(1, 0, 0)$ & \\
      \hline
     Nonlinear separatrix & 4 & $(1, 0, 0, 0, 0)$ & $(1, 0, 0, 0, 0)$ & $(1, 0, 0, 0, 0)$ & & \\
      \hline
      Hill system with PO & 3 & $(1, 0, 0, 0)$ & $(1, 0, 0, 0)$ & $(1, 0, 0, 0)$ & $(1, 0, 0, 0)$ & $(1, 0, 0, 0)$ \\
      \hline 
     EMT Hill system& 6 &  $(1, 0, 0, 0, 0, 0, 0)$ & $(1, 0, 0, 0, 0, 0, 0)$ & $(1, 0, 0, 0, 0, 0, 0)$ & & \\
      \hline
     Ellipsoidal bistable & 2 & $(1, 1, 0)$ & $(1, 1, 0)$ & $(1, 0, 0)$ & & \\
      \hline
     Ellipsoidal bistable& 3 & $(1, 0, 1, 0)$ & $(1, 0, 1, 0)$ & $(1, 0, 0, 0)$ & & \\
      \hline
     Ellipsoidal bistable& 4 & $(1, 0, 0, 1, 0)$ & $(1, 0, 0, 1, 0)$ & $(1, 0, 0, 0, 0)$ & & \\
      \hline
     Ellipsoidal bistable& 5 & $(1, 0, 0, 0, 1, 0)$ & $(1, 0, 0, 0, 1, 0)$ & $(1, 0, 0, 0, 0, 0)$ & & \\
     \hline
    \end{tabular}
    \caption{The correct Betti numbers for the sets $U$ and $N_k$ for each example.}
    \label{tab:correct-homology}
\end{table}

\subsection{Hill system with a periodic orbit}
\label{sec:p3d}
We define a dynamical system on 
$X=[0,3.062] \times [0,4.072] \times [0,11.26362]$.
Let $n = 10$, $\gamma = [1, 1, 1]^T$,

\[ 
L=
\begin{bmatrix}
0.133 & 1.743 & 0\\
0.599 & 0.124 & 0.245\\
0.175 & 0 & 0.458
\end{bmatrix},
U=
\begin{bmatrix}
0.578 & 3.299 & 0\\
1.994 & 0.725 & 5.205\\
1.737 & 0 & 2.164
\end{bmatrix},
~\text{and}
\]

\[
\Theta = 
\begin{bmatrix}
0.363 & 1.531 & 0\\
2.036 & 0.862 & 0.233\\
0.818 & 0 & 0.168
\end{bmatrix}.
\]

Define
\begin{equation}
\makeatletter
\scalebox{0.9}{%
$\displaystyle
\begin{cases}
h_{1}(x) = H^{+}(x_1, L_{11}, U_{11}, \Theta_{11}, n) H^{-}(x_2, L_{21}, U_{21}, \Theta_{21}, n) H^{-}(x_3, L_{31}, U_{31}, \Theta_{31}, n) \\
h_2(x) = H^{+}(x_{2}, L_{22}, U_{22}, \Theta_{22}, n) H^{-}(x_1, L_{12}, U_{12}, \Theta_{12}, n) \\
h_3(x) = H^{+}(x_{3}, L_{33}, U_{33}, \Theta_{33}, n) H^{-}(x_{2}, L_{23}, U_{23}, \Theta_{23}, n)
\end{cases}
$%
}
\makeatother
\end{equation}
where
\begin{equation}
\label{eq:neg_Hill}
H^{-} (y, \ell, u, \theta) =  \ell + (u - \ell) \frac{\theta^n} {y^n + \theta^n},
\end{equation}
and
\begin{equation}
\label{eq:pos_Hill}
H^{+} (y, \ell, u, \theta) = \ell + (u - \ell) \frac {y^n} {y^n + \theta^n}.
\end{equation}

Define the system $\dot{x} = f(x)$ componentwise by
\begin{equation}
\label{eqn:periodic_3d}
\dot{x_i} = - \gamma_i x_i + h_i(x), \quad i = 1, 2, 3.
\end{equation}

\subsection{Six-dimensional EMT Hill system}

We now define a dynamical system on $X=[0,0.5] \times [0,3.558] \times [0,1.57345566] \times [0,5.762] \times [0,15.169196784] \times [0,5.056847]$. 
Let $n = 10$, $\gamma = [1, 1, 1, 1, 1, 1]^T$, 
\[  
L=
\begin{bmatrix}
0 & 0 & 0.318 & 0 & 0 & 0\\
0.148 & 0 & 0 & 0 & 0.610 & 0\\
0 & 1.210 & 0.510 & 0 & 1.054 & 0.993\\
0.150 & 0 & 0 & 0 & 0.6454 & 0\\
0 & 1.024 & 0 & 0.707 & 0 & 0.684\\
0 & 0 & 0.533 & 0 & 0 & 0
\end{bmatrix},
\]
\[ 
U=
\begin{bmatrix}
0 & 0 & 1.316 & 0 & 0 & 0\\
0.520 & 0 & 0 & 0 & 2.037 & 0\\
0 & 1.442 & 0.965 & 0 & 2.056 & 3.917\\
0.670 & 0 & 0 & 0 & 3.622 & 0\\
0 & 2.113 & 0 & 3.720 & 0 & 1.291\\
0 & 0 & 1.239 & 0 & 0 & 0
\end{bmatrix},
\]
and
\[ 
\Theta =
\begin{bmatrix}
0 & 0 & 0.250 & 0 & 0 & 0\\
1.312 & 0 & 0 & 0 & 1.779 & 0\\
0 & 0.668 & 0.206 & 0 & 0.478 & 0.166\\
1.333 & 0 & 0 & 0 & 2.881 & 0\\
0 & 5.698 & 0 & 1.676 & 0 & 0.950\\
0 & 0 & 0.705 & 0 & 0 & 0
\end{bmatrix}.
\]
        
\begin{equation}
\makeatletter
\scalebox{0.9}{%
$\displaystyle
\begin{cases}
h_1(x) = H^{-}(x_2, L_{21}, U_{21}, \Theta_{21}, n) H^{-}(x_4, L_{41}, U_{41}, \Theta_{41}, n) \\
h_2(x) = H^{-}(x_{3}, L_{32}, U_{32}, \Theta_{32}, n) H^{-}(x_5, L_{52}, U_{52}, \Theta_{52}, n) \\
h_3(x) = H^{+}(x_1, L_{13}, U_{13}, \Theta_{13}, n) H^{-}(x_3, L_{33}, U_{33}, \Theta_{33}, n) H^{-}(x_6, L_{63}, U_{63}, \Theta_{63}, n) \\
h_4(x) = H^{-}(x_5, L_{54}, U_{54}, \Theta_{54}, n) \\
h_5(x) = H^{+}(x_3, L_{35}, U_{35}, \Theta_{35}, n) H^{-}(x_2, L_{25}, U_{25}, \Theta_{25}, n) H^{-}(x_4, L_{45}, U_{45}, \Theta_{45}, n) \\
h_6(x) = H^{-}(x_{3}, L_{36}, U_{36}, \Theta_{36}, n) H^{-}(x_{5}, L_{56}, U_{56}, \Theta_{56}, n),
\end{cases}    
$%
}
\makeatother
\end{equation}
where $H^{\pm}$ are defined by \eqref{eq:neg_Hill} and \eqref{eq:pos_Hill}.

Define the system $\dot{x} = f(x)$ componentwise by
\begin{equation}
\label{eqn:EMT}
\dot{x_i} = - \gamma_i x_i + h_i(x), \quad i = 1, 2, \ldots, 6.
\end{equation}

\end{document}